\def\aff           {affine Lie algebra}
\def\alg           {algebra}
\def\apc           {{{\mathfrak A}_{\tC,\vec p}}}
\def\apcg          {{{\mathfrak A}_{\tC_g,\vec p}}}
\def\auto          {automorphism}
\newcommand\BArray[4]{\mbox{\large$[$}\!\!{\scs\begin{array}{lr}{}\\[-1.94em]
                   {\scs #1}\!\!&\!\!{\scs #2}\\[-.43em]{\scs #3}\!\!%
                   &\!\!{\scs #4}\\[-.4em] \end{array}}\!\!\mbox{\large$]$}}
\def\bc            {boundary condition}
\def\Bc            {Boundary condition}
\def\be            {\begin{equation}}
\def\bearl         {\begin{array}{l}}
\def\bearll        {\begin{array}{ll}}
\def\bearlll       {\begin{array}{lll}}

\def\betab         {{\bar\beta}}
\def\bfe           {{\bf1}}

\def\bub           {{\bar\beta}}
\def\calb          {{\cal B}}
\def\calbg         {{\cal B}^{\sss(g)}}
\def\calc          {{\cal C}}
\def\calf          {{\cal F}}
\def\calg          {{\cal G}}
\def\calh          {{\cal H}}
\def\calhb         {\bar{\cal H}}
\def\calhmv        {\vec{\cal H}_\muv}

\def\calm          {{\cal M}}
\def\calo          {{\cal O}}
\def\cals          {{\cal S}}
\def\calu          {{\cal U}}
\def\cb            {chiral block}
\def\cft           {conformal field theory}
\def\Cft           {CFT}
\def\cfts          {conformal field theories}
\def\chie          {\chii^{\sss(0)}}
\def\chii          {\raisebox{.15em}{$\chi$}}
\def\chiO          {\chii^{{\cal O}}}
\def\chir          {\mbox{$\mathfrak A$}}
\def\chiR          {{\mathfrak A}}
\def\chirb         {\mbox{$\bar{\mathfrak A}$}}
\def\chirB         {\bar{\mathfrak A}}
\def\chirg         {\mbox{${\mathfrak A}^G$}}
\def\chiz          {\chii^{\sss(1)}}
\def\cla           {classifying algebra}
\def\clA           {\mbox{$\calc(\chiR,\chiR)$}}
\def\clAb          {\mbox{$\calc(\chiR,\chirB)$}}
\def\class         {classification}

\def\complex       {{\mathbb C}}
\def\con           {conformal }
\def\Con           {Conformal }
\def\corfu         {correlation function}

\def\dyd           {Dynkin diagram}
\def\eal           {E^{\bar\alpha}}
\def\ebe           {E^{\betab}}
\def\eE            {{\rm e}}
\def\eear          {\end{array}}

\def\eq            {\,{=}\,}
\newcommand\erf[2] {(\ref{#1#2})}
\newcommand\F[6]   {\mbox{\sf F}_{\! #1,#2}\BArray{#4}{#5}{#3}{#6}}
\def\findim        {finite-dimensional}
\def\fiss          {\varphi_{1,s}}
\newcommand\Frac[2]{\mbox{\large$\frac{#1}{#2}$}}
\def\fpc           {{{\cal O}(\tC{\setminus}\vec p)}}
\def\fpe           {{{\cal O}(\Pe{\setminus}\vec p)}}
\def\ft            {field theory}
\def\furu          {fusion rule}
\def\g             {\mbox{$\liefont g$}}
\def\G             {{\rm G}}
\def\gb            {\mbox{$\bar{\liefont g}$}}

\def\gs            {{\sigma^*}}
\def\Gs            {{G^*_{\phantom|}}}
\def\Hat           {\hat }
\def\hopf          {Hopf algebra}
\newcommand\hsp[1] {\mbox{\hspace{#1 em}}}
\def\hy            {$\mbox{-\hspace{-.66 mm}-}$}
\def\id            {\mbox{\sl id}}

\def\ii            {{\rm i}}
\def\iN            {\,{\in}\,}
\def\infdim        {infinite-dimensional}
\def\Infdim        {Infinite-dimensional}
\def\intro         {introduction }
\def\Intro         {Introduction }
\def\irrep         {irreducible representation}

\def\J             {{\rm J}}
\def\Jp            {{\J^{+}_{\phantom i}}}
\def\K             {{\rm K}}
\def\kma           {Kac\hy Moody algebra}
\def\kzc           {Knizh\-nik\hy Za\-mo\-lod\-chi\-kov connection}
\def\L             {\mbox{$\mathfrak L$}}
\newcommand\labl[2]{\label{#1#2}}

\def\lambdad       {{\dot\lambda}}
\def\lambdap       {{\lambda^{\!+}_{\phantom i}}}

\def\lie           {Lie algebra}
\def\Lie           {Lie group}

\def\liefont       {\mathfrak }
\def\llb           {\mbox{\large(}}
\def\Llb           {\mbox{\Large[}}

\def\lrb           {\mbox{\large)}}
\def\Lrb           {\mbox{\Large]}}

\def\Ltimes        {{\mbox{\normalsize$\ltimes$}}}
\def\Lv            {L\raisebox{.52em}{$\sss\!\vee$}}
\def\Lw            {L\raisebox{.52em}{${\sss\!\vee}\scs*$}}
\def\lwv           {L_{\rm w}\raisebox{.51em}{$\!\!\!\!\sss\vee$\,}}

\def\modinv        {modular invarian}
\def\Modinv        {Modular invarian}
\def\mU            {{[\mu]}}
\def\Mu            {{[\mu,\psu]}}
\def\mub           {{\bar\mu}}
\def\mud           {{\dot\mu}}
\def\mup           {{\mu^{\!+}_{\phantom i}}}
\def\muP           {{\mu^{\!+}_{\phantom.}}}
\def\muv           {{\vec\mu}}
\def\MV            {\dot L\raisebox{.51em}{$\sss\!\vee$}}
\def\Mw            {\dot L\raisebox{.51em}{${\sss\!\vee}\scs*$}}
\newcommand\N[3]   {{\rm N}_{#1,#2}^{\ \ \ #3}}
\def\nE            {\,{\ne}\,}

\newcommand\nxt[1] {\\\raisebox{.12em}{\rule{.35em}{.35em}}\hsp{.6}#1}
\newcommand\Nxt[1] {\\[.1em]\raisebox{.12em}{\rule{.35em}{.35em}}\hsp{.6}#1}
\newcommand\nxT[1] {\raisebox{.12em}{\rule{.35em}{.35em}}\hsp{.6}#1}

\def\nub           {{\bar\nu}}
\def\nup           {{\nu^{\!+}_{\phantom i}}}
\def\nuP           {{\nu^{\!+}_{\phantom,}}}
\newcommand\ol[1]  {^{[#1]}}
\def\omego         {\sigma^{}_{\!\circ}}
\def\Omego         {\sigma_{\!\circ}}
\def\ommuv         {\sigma^{}_{\vec{\bar\nu}}}
\def\Ommuv         {\vec\sigma^{}_{\vec{\bar\nu}}}
\def\Ommuvs        {\vec\sigma^*_{\!\vec{\bar\nu}}}
\def\one           {\mbox{\small $1\!\!$}1}
\def\onedim        {one-dimen\-sional}
\def\onetom        {1,2,...\,,m}
\def\oplu          {\,{\oplus}\,}
\def\ot            {\raisebox{.07em}{$\scriptstyle\otimes$}}
\def\oT            {\,\ot\,}

\def\parfu         {partition function}
\def\Pe            {\mathbb{CP}^1}
\newcommand\pho[1] {\phi^{}_{#1,\tilde{#1}}}
\newcommand\Pho[1] {\phi^{}_{#1,{#1}^+_{\phantom .}}}
\def\pp            {p}
\def\pp            {p\!\!p}
\def\pslz          {\mbox{PSL(2,$\zet$)}}
\def\psu           {{\psi}}
\def\q             {quantum }
\def\Q             {Quantum }
\def\qft           {quantum field theory}
\newcommand\rca[3] {{\rm R}^{#1_{\phantom e}}_{#2,\tilde{#2}}}
\newcommand\rch[3] {{\rm R}^{#1_{\phantom.}}_{\Hat #2}}

\def\reals         {{\mathbb R}}
\def\rep           {representation}
\def\Rep           {Representation}
\def\Res           {{\cal R}\!\mbox{\sl es}}
\def\resp          {respectively}
\def\rhob          {{\bar\rho}}

\newcommand\sect[1]{\section{#1}\setcounter{equation}{0}}
\def\scs           {\scriptstyle}

\def\Se            {S^{\sss(0)}}
\def\sigmas        {\sigma_{\!\ss}}
\def\SJ            {S^\J}
\def\slie          {Lie su\-per\-al\-ge\-bra}
\def\slz           {\mbox{SL(2,$\zet$)}}
\def\so            {\mathfrak{so}}
\def\SO            {S^\calo}
\def\sss           {\scriptscriptstyle}
\def\ss            {{\rm s}}

\def\sus           {{I}}
\def\sush          {{\Hat\sus}}

\def\sym           {sym\-me\-try}
\def\syms          {sym\-me\-tries}
\def\Sz            {S^{\sss(1)}}
\def\tC            {{\hat X}}
\def\Tau           {\Theta}
\def\Te            {T^{\sss(0)}}
\def\tft           {topological field theory}

\def\TN            {\tilde{\rm N}}
\newcommand\tNC[3] {\Hat{\rm N}_{#1,#2}^{\ \ \ #3}}
\newcommand\tNl[3] {\Hat{\rm N}_{\Hat #1,\Hat #2,\Hat #3}}
\newcommand\TNl[6] {\Hat{\rm N}_{(#1,#2),(#3,#4),(#5,#6)}}
\newcommand\tNN[3] {\Hat{\rm N}_{\Hat #1,\Hat #2}^{\ \ \ \Hat #3}}

\def\TO            {T^\calo}
\def\twodim        {two-dimensional}
\def\tS            {\hat S}

\def\Tz            {T^{\sss(1)}}

\def\vac           {\Omega}
\def\vacb          {{\bar\Omega}}
\def\Vee           {{\sss\vee}}

\def\voa           {vertex operator algebra}
\def\vop           {vertex operator}
\def\vv            {{\rm v}}
\def\vvir          {v_{{\rm vir}}}
\def\wrt           {with respect to }
\def\wrtt          {with respect to the }
\def\wzwm          {WZW model}

\def\xxt           {^{{\sss\rm xt}_{\!\!\!\phantom:}}}
\def\Xxt           {\raisebox{.59em}{$\sss\rm xt$}}
\def\XXt           {\raisebox{.59em}{$\sss\rm xt\!\!\!\!\!\!$}}
\def\zet           {{\mathbb Z}}

\def\zetpluso      {{\mathbb Z}_{\ge0}} 

\documentclass[12pt]{article} \usepackage{amssymb,amsfonts,latexsym}

\begin{document}
\begin{flushright}  {~} \\[-1cm]
{\sf math.QA/0011160}\\[1mm]{\sf PAR-LPTHE 00-40}\\[1mm]
{\sf November 2000} \end{flushright}
 
\begin{center} \vskip 14mm
{\Large\bf LIE ALGEBRA AUTOMORPHISMS}\\[3mm]
{\Large\bf IN CONFORMAL FIELD THEORY}
\raisebox{.66em}{\rule{.55em}{.55em}} \\[19mm]
{\large J\"urgen Fuchs$\;^1$ \ and \ Christoph Schweigert$\;^2$}
\\[8mm]
$^1\;$ Institutionen f\"or fysik~~~~{}\\
Universitetsgatan 1\\ S\,--\,651\,88\, Karlstad\\[5mm]
$^2\;$ LPTHE, Universit\'e Paris VI~~~{}\\
4 place Jussieu\\ F\,--\,75\,252\, Paris\, Cedex 05
\end{center}
\vskip 18mm
\begin{quote}{\bf Abstract}\\[1mm]
The role of automorphisms of infinite-dimensional Lie algebras in
conformal field theory is examined. Two main types of applications are
discussed; they are related to the enhancement and reduction of 
symmetry, respectively. The structures one encounters also
appear in other areas of physics and mathematics. In particular, they 
lead to two conjectures on the sub-bundle structure of chiral blocks,
and they are instrumental in the study of conformally invariant 
boundary conditions.
\end{quote}
\vfill
$\!\!\!\!\!\raisebox{.32em}{\rule{.35em}{.35em}}$ \
Invited talk by J.\ Fuchs at the conference on {\sl Infinite-Dimensional
Lie Theory and Conformal\\
Field Theory\/} (University of Virginia,
May 2000, {\tt http:/$\!$/www.math.virginia.edu/ilacft/}).\\
To appear in Con\-temp.\,Math.
 

\sect{Observables and \auto s in \qft}

In quantum theory, a physical system is described in terms of its space 
$\calh$ of physical states -- some vector space over $\complex$\,, in many
situations a Hilbert space -- and of {\em fields\/} acting as operators on
$\calh$. Distinguished among the fields are the {\em observables\/}.
In the quantum mechanical description of a point particle, for example,
$\calh$ is a Hilbert space and the observables are densely defined
operators on $\calh$. One may wish to endow the collection of all observables 
with the structure of a unital associative \alg. But since two operators 
can only be multiplied if the image of the first operator is contained in 
the range of the second, this idea is somewhat too naive already in quantum 
mechanics. 

Once one studies quantum field theory rather than quantum mechanics, the 
algebraic structure must be further refined. For instance, there is a 
framework \cite{HAag} in which the local observables of relativistic \qft\ 
give rise to a net of von Neumann \alg s. Tentatively, the space $\calh$ 
of physical states may then be regarded as a module over this (net of) 
algebra(s), which may decompose into a direct sum $\bigoplus_\mu\!\calh_\mu$
of irreducible modules $\calh_\mu$.
But a refinement of this point of view proves to be necessary, involving the
concept of {\em superselection rules\/}. This is briefly expressed by
saying that vectors in different $\calh_\mu$ must not be linearly combined;
a precise formulation is that one must not regard the space of states
as the direct sum of modules $\calh_\mu$ over
the observable algebra, but rather just as a collection of those modules.
The spaces $\calh_\mu$ and the associated \rep s $R_\mu$ of the observables
are then called superselection sectors \cite{STwi}. Thus
observables act as operators within the individual modules $\calh_\mu$, 
whereas non-observable (`charged') fields operate between different modules.
(An example of such a field operator is provided by the fermion field $\psi$, 
which changes the fermion number and thereby the superselection sector.)

The spaces $\calh_\mu$ are not just modules over the observable \alg, but 
enjoy additional properties. For instance, to incorporate
compatibility with the probabilistic features of quantum physics,
they should be unitarizable. Also, one would like to have a sufficiently 
large observable \alg, and hence as a kind of naturality property 
requires the modules $\calh_\mu$ to be simple. On the other hand, it is 
economic if the observable \alg\ is not too large either; this is 
implemented by a {\em sta\-te-field correspondence\/} according 
to which any vector in each of the spaces $\calh_\mu$ is obtainable by 
acting with a unique field operator on a unique (up to a scalar) {\em 
vacuum vector\/} $v_\vac$. The module for which $v_\vac$ is a cyclic vector 
is called the vacuum superselection sector and is denoted by $\calh_\vac$.
Now for any endomorphism $\varpi$ of the observable \alg\ and any \rep\ $R$,
the composition $R\,{\circ}\,\varpi$ is again a \rep\, on the same vector 
space; it inherits crucial features from $R$, such as simplicity and, for 
suitable $\varpi$, unitarizability. Whether the \rep\ $R\,{\circ}\,\varpi$ 
is physically interesting depends on the endomorphism $\varpi$ (and on $R$, 
of course).  In \qft, the issues of having an interesting \rep\ and of having 
an interesting endomorphism are intimately related. Indeed, there should exist
a collection of special endomorphisms $\varpi_\mu$ that satisfy
$R_\mu\,{\cong}\,R_\vac{\circ}\,\varpi_\mu$, with $R_\vac$ the
{\em vacuum \rep\/} that corresponds to $\calh_\vac$.
All these structures can be made explicit in the formalization of local
observables by nets of von Neumann algebras alluded to above.

In that framework, via composition of the associated endomorphisms one can also 
endow the collection of \rep s $R_\mu$ with a tensor product, which is fully 
reducible. On the set of isomorphism classes of simple modules $\calh_\mu$ 
-- to be referred as {\em sectors\/}, for short, and to be denoted by 
$\lambda,\,\mu,...$ -- this tensor product provides a multiplication called the 
fusion product and denoted by `$\star$'. By full reducibility one can write
  \be  \lambda\star\mu = \sum_{\nu\in\sus} \N\lambda\mu\nu\,\nu \,,
  \labl1N  \end{equation}
where $\sus$ stands for the set of sectors.
The vacuum sector $\vac$ is a unit element for the fusion product.
Also, the evaluation on the unit element, $C_{\lambda,\mu}
\,{:=}\,\N\lambda\mu{\!\vac}$, provides an involutive \auto\ of the fusion 
rules (otherwise the QFT model wouldn't possess sensible two-point \corfu s). 
We write $C{:}\;\mu\,{\mapsto}\,\mup$, i.e.\ 
$C_{\lambda,\mu}^{}\eq\delta_{\lambda,\mup}$. In short, one has the
structure of a {\em fusion ring\/} -- a commutative associative unital ring 
$\cal R$ over $\zet$, with a distinguished basis in which the structure constants 
are non-negative integers $\N\lambda\mu\nu$, and with evaluation at the unit 
element being a conjugation. (Instead of $\cal R$, one may also consider
${\cal R}\,{\otimes}^{}_\zet\complex\,$, called the {\em fusion algebra\/}.
The relations \erf1N are also referred to as the {\em fusion rules\/}.)

Of special interest are {\em rational\/} theories, which 
possess only a finite number of sectors, $|\sus|\,{<}\,\infty$.
A feature special to the rational case is the existence of a
matrix $S$ that diagonalizes the fusion rules in the sense that
  \be  \N\lambda\mu\nu = \sum_{\kappa\in I} {S_{\kappa,\lambda}^{}\,
  S_{\kappa,\mu}^{}\,S_{\kappa,\nu}^*} \,/\, {S_{\kappa,\vac}} \,.  \labl1v
  \end{equation}

The two \rep s $R$ and $R\,{\circ}\,\varpi$ are particularly close 
relatives when $\varpi\eq\omega$ is an {\em auto\/}morphism. 
Sectors associated via $R_\mu\,{\cong}\,R_\vac{\circ}\,\omega_\mu$ to \auto s 
$\omega_\mu$ are therefore distinguished; they are called {\em simple currents\/}.
Trivially, there is always at least one simple current, the vacuum sector $\vac$.
The simple currents $\J$ are the elements of the distinguished basis of the
fusion ring $\cal R$ that are units of $\cal R$. They are 
characterized by $\sum_\mu\!\N\lambda\J\mu\eq1$ for all $\lambda\iN\sus$; 
under the fusion product they form a subring of $\cal R$ that is the
group ring of an abelian group. In rational theories,
simple currents can also be characterized by $S_{\J,\vac}\eq S_{\vac,\vac}$.

\medskip

Let us mention that the observables can be regarded as encoding part of the
{\em symmetry\/} of a physical system. (They play in particular the role of
a `spectrum generating symmetry'.)
Now `symmetry' is one of the most helpful concepts in physics,
but it shares the fate of other physical 
notions in that by no means does it refer to a single well-identified 
concept, but rather encompasses a whole zoo of different ideas.
E.g.\ a common incarnation of symmetry that is not covered by the above
discussion of observables arises in the form of transformations that act 
non-trivially on the field variables used to describe a physical system, but 
leave the system itself invariant. 
Concrete examples include various groups of transformations of the
canonical variables in Lagrangian mechanics and field theory.\,%
 \footnote{~The set of all such transformations then has the structure of a 
 group. In quantum theory a Lagrangian description of a system need not exist,
 though. Accordingly, more general structures than groups can occur as symmetry
 transformations, see e.g.\
 \cite{frrs,rehr2,masc7,fugv3,bons}. Collectively these structures are often
 referred to as `quantum symmetries' or (in a non-technical sense) quantum groups.
 While these issues do not concern us here, in any case we use the term `\qft'
 in a broad sense,
 not implying a Lagrangian setting, even not the existence of a classical limit
 (and hence in particular not restricted to the realm of perturbation theory).}
   
\sect{\lie s and their \auto s in \cft}

We now focus our attention on a particular class of models -- rational models 
of {\em \twodim\ \cft\/} ({\em CFT\/}). The conformal symmetry (and extensions 
thereof) present in such models is encoded in a suitable algebra of observables. 
The structure of CFT is severely constrained by the geometry of the underlying
{\em world sheet\/}, i.e.\ the two-dimensional manifold $X$ on which the 
field theory is defined. In particular, 
the observables split into quantities depending either analytically or
anti-analytically on a local complex coordinate on the world sheet; these are
referred to as {\em chiral\/} and {\em anti-chiral\/} observables, \resp.
Studies of special classes of models (see e.g.\ \cite{witt39}) indicate 
that a similar splitting happens for other quantities such as \corfu s; this 
phenomenon is known as `holomorphic factorization'.
In free field models, and consequently also in the string theories constructed
from them, this amounts to a decomposition into `left-' and `right-movers'.

Because of their interest in string theory, we consider manifolds $X$ that
are compact surfaces of Euclidean 
signature. In order to constitute a valid world sheet for CFT, such a surface 
$X$ must be endowed with a conformal structure, i.e.\ an equivalence class, 
\wrtt local rescalings, of metrics. $X$ may have a boundary or be 
non-orientable. In the orientable case, together with a choice of orientation 
every conformal structure provides a complex structure. The best way to 
account for the splitting into chiral and anti-chiral observables on 
arbitrary world sheets is to replace \cite{scfw} the physical world sheet $X$ 
by its Schottky double cover $\tC$, which is a complex curve.  The theory 
to be studied on $\tC$ is called {\em chiral CFT\/}, or {\em CCFT\/}
for short; until the end of section 5, we will mainly deal with CCFT.

There exist several different formalizations of the notion of chiral CFT, 
which are not entirely equivalent. Here we adopt an approach that is 
well suited for various applications to string theory as well as to
statistical mechanics and condensed matter physics.\,%
 \footnote{~Other approaches can e.g.\ be found in \cite{frrs,EVka,rehrX} 
 (local observables forming nets of von Neumann algebras over the circle) and 
 in \cite{gago}.
 For a functional-analytic interpretation of VOAs see \cite{huan10}.}
Its basic features are the following:
\nxt
The observables of a CCFT model have the structure of a rational
{\em vertex operator \alg\/} ({\em VOA}), or conformal vertex \alg, \chir.
Such algebras, in turn, can for instance be constructed \cite{KAc4} from 
so-called Lie algebras \L\ of formal distributions.
\nxt
Along with the VOA \chir\ comes its \rep\ category Rep(\chir), which is 
a se\-mi-simple tensor category.
\nxt
The isomorphism classes of simple \chir-modules are the sectors of the model.
\nxt
There is a lot more structure on the \rep\ category Rep(\chir). It is 
conjectured\,%
 \footnote{~It is so far not known whether the \rep\ category of every rational
 VOA indeed possesses all properties of a modular tensor category, though this
 property has been established for several classes of VOAs, compare e.g.\
 \cite{huan7,hule5}. In the approach based on nets of von Neumann algebras, the
 presence of a (unitary) modular tensor category has been established
 in \cite{kalm}.}
that this can be encoded in so-called Moore\hy Seiberg data \cite{mose3,fefk3,bakI}, 
which equip Rep(\chir) with the structure of a {\em modular tensor category\/}. 
\nxt
The presence of a modular tensor category is, in turn, the basis for the
relation of CFT to three-dimensional topological field theory (TFT),
see e.g.\ \cite{frki,mose4,witt27,fffs3}, and thereby
to invariants of knots and three-manifolds \cite{TUra}.  

Note that this formalization does not entirely fit into
the general field theoretic framework of section 1; for example we are
not given a relation between the sectors and endomorphisms of the observables.
Nevertheless the important features mentioned there are again present. For 
instance, the vertex operator map supplies a state-field correspondence in the 
vacuum sector. Also, a modular tensor category is in particular spherical, so 
it has a Grothendieck ring, and this ring precisely realizes the fusion rules.  
The tensor identity object is the vacuum sector $\vac$, for which 
the underlying vector space $\calh_\vac$ is the vector space of the VOA 
($\calh_\vac$ naturally carries the structure of a module over the VOA itself).
For more details we refer to \cite{freN5,BAki,scfw} and the literature cited there. 

As we are concerned with rational theories, $|\sus|\,{<}\,\infty$,\,%
 \footnote{~However, many of the concepts and results for rational CFT in fact
 carry over to large classes of non-rational ones.}
there is a matrix $S$ that diagonalizes the fusion rules, see
\erf1v. The requirement that the CCFT gives rise to a modular tensor category
includes the statement that this matrix $S$ and the diagonal matrix $T$ -- that 
has entries $T_{\mu,\nu}\eq\delta_{\mu,\nu}T_\mu$, $T_\mu\eq\exp(2\pi\ii(\Delta
_\mu\,{-}\,c/24))$, with $\Delta_\mu$ the conformal weights of the sectors and
$c$ the eigenvalue of the Virasoro central charge -- satisfy the defining 
relations $S^4\eq\one$ and $(ST)^3\eq S^2$ for the generators of \slz.  

\smallskip 

Next we comment on the \infdim\ complex \lie s \L\ of formal distributions 
from which one can construct \voa s. The idea in that construction is to 
define a VOA \chir\ such that \L\ can be regarded as the \lie\ of Laurent 
modes for the vertex operators associated to specific vectors in \chir.
The most prominent examples for such \lie s of formal distributions are the
{\em Virasoro \alg\/} -- the relevant vector in \chir\ is the Virasoro element
$\vvir$ --\,, the Heisenberg algebra, and {\em untwisted affine \lie s\/}. 
Various other classes of such \lie s \L\ are known. Interesting
results about the corresponding CFT models can often be derived by relating 
them to simpler models that are based on the \alg s just mentioned, by
special constructions like the `coset construction' or `Hamiltonian reduction'.
But the chiral algebras for these models are 
less understood and their \rep\ theory is not fully developed.  

One should think of \chir\ \resp\ \L\ as expressing the chiral symmetries 
of a CFT model {\em locally\/} in a neighborhood of any point on the
curve $\tC$. One also wants to express the symmetries {\em globally\/} on 
$\tC$; this is achieved by imposing appropriate invariance conditions, known 
as {\em Ward identities\/}, on the \cb s.  When \chir\ is obtained from a 
\lie\ \L\ of formal distributions, then the Ward 
identities frequently reduce to a (co-)invariance condition \wrt another type
of \infdim\ \alg s $\apc$, which we like to call {\em block \alg s\/}\,%
 \footnote{\,In \cite{beDr3}, the term `chiral \alg' is used for the objects
 $\apc$. This use of the term is unfortunately different from terminology in
 most of the physics literature, where it is the VOA \chir\ itself or the Lie
 \alg\ \L\ of Laurent modes that is called chiral \alg.} 
(for more details, see section 4).
The block \alg s can be described in terms of globally defined quantities like 
functions and differential forms on $\tC$, or sections in corresponding bundles.
The \alg s \chir\ and \L\ are related to them via Laurent expansion of such 
global quantities in a local holomorphic coordinate.

An important class of CFT models are the {\em WZW models\/}. They provide the 
input for the coset construction and for Hamiltonian reduction. WZW models 
possess a Lagrangian field theory realization, as sigma models with target 
space a (connected simply connected compact) real Lie group \G. The \voa\ of 
a WZW model is generated by a \lie\ of formal 
distributions, which is just an untwisted affine \lie. We denote this affine 
\lie\ by \g; the horizontal sub\alg\ \gb\ of \g\ is the complexification 
of the \lie\ ${\rm Lie}(\G)$ of the group manifold \G. The block \alg s of 
WZW models are again \infdim\ \lie s, obtained as tensor products of \gb\ 
with certain \alg s\ of meromorphic functions on $\tC$, see formula \erf gf.

\smallskip 

Below we will examine two typical situations in which \auto s enter in \qft:
  \begin{enumerate}\addtolength\itemsep{.04em}
 \item The enhancement of symmetries.
 \item The restriction of observables.
  \end{enumerate}
\vskip.2em
In CFT, the first issue
arises in the form of extended chiral algebras, notably so-called 
{\em simple current extensions\/}, and the second in the form of {\em orbifold\/}
constructions. The operations of forming the orbifold of a CFT \wrt a 
finite abelian group and of simple current extension are actually each other's 
inverse. (Non-abelian orbifolds are inverse to more general extensions.)
Thus in principle there is no need to study the two constructions separately. 
Nevertheless in practice it is most sensible to do so.
The reason is that some CFT models -- notably WZW models -- are much better 
understood than others and hence serve as natural starting points for each
of the two constructions.

These matters will be studied in sections 3 and 5, respectively. Let us include 
into the list as separate items also two more specific contexts -- to be dealt
with in sections 4 and 7 -- in which these structures are encountered:
  \begin{enumerate}\addtolength\itemsep{.04em}
 \item[3.] The sub-bundle structure of chiral blocks.
 \item[4.] Symmetry reduction caused by the presence of boundaries.
  \end{enumerate} 
\vskip.2em
The role of \auto s in these situations is the following.
  \begin{enumerate}\addtolength\itemsep{.32em}
 \item Symmetry enhancement: Simple currents naturally correspond to
  (classes of) outer \auto s of the \lie\ \L\ of formal distributions.
  (They do not preserve the Virasoro algebra and hence do not 
  give rise to \auto s of the VOA \chir.) \\
  For instance, in the WZW case, where $\L\eq\g$ is an untwisted affine \lie,
  it can be proven \cite{jf15,fugv2} that (barring a single exceptional case
  present for $E_8$ level 2) there is a natural isomorphism
    \be  \calg(\g) \cong {\cal Z}(\G)  \labl gh \end{equation}
  between the simple current group $\calg$ and
  the center of the relevant simply connected Lie group \G. 
 \item Restriction of observables: The orbifold group is 
  a group of \auto s of \L. The orbifold model should still possess
  conformal symmetry; accordingly, the fixed point set of the orbifold
  group must contain the Virasoro subalgebra of \L.
  As a consequence, in this case the \auto s of \L\ induce \auto s of \chir.
 \item Chiral blocks: Suitable outer \auto s of \L\ give rise to finite order
  \auto s of the bundles of chiral blocks. These bundles then split into
  sub-bundles invariant under those \auto s.
 \item Symmetry breaking \bc s: In CFT they are partially characterized by the
  subalgebra of \chir\ \resp\ \L\ that they preserve. (Again for preserving 
  conformal invariance this must contain the Virasoro \alg.) Their analysis 
  is simplified technically when that subalgebra can be characterized as 
  the fixed point subalgebra under some group of \auto s.
  \end{enumerate}

To summarize these introductory remarks, let us repeat that in \Cft\ (and in 
other areas of physics as well) various structures related to \auto s of \lie s 
appear naturally. Which \auto s are of interest depends largely on the concrete 
system under study.  In this paper we restrict our attention to the 
applications listed above. (The exposition mainly reviews the subject, but 
a few results, such as formula \erf39, have not appeared previously.) 
Before proceeding to those applications, we also should like to mention
that in many circumstances, in particular when it comes to model-dependent
issues, in the study of \Cft\ the arguments are not, or not yet, 
rigorous, but involve some amount of heuristics -- often concealed as
{\em physical intuition\/}. Nevertheless, there do exist situations in which 
one can make precise statements and establish rigorous proofs. A well-known 
illustration is the WZW Verlinde formula, which will be mentioned in section 4. 
Below we will indicate the relevant literature where such proofs are 
available, but not go into the details of any of the proofs.

\smallskip

When studying the \auto s of our interest
we will also encounter the following general structures.
First, given an \L-\rep\ $R$ on the vector space $V$ and an \auto\ $\omega$ of \L, 
there is a `twisted intertwiner', i.e.\ a linear map $\Tau_\omega\,{\equiv}\,\Tau_
\omega^{(R)}$ from the \L-module $(R,V)$ to $(R{\circ}\omega,V)$ satisfying
  \be  (R{\circ}\omega)(x) \circ\Tau_\omega = \Tau_\omega\circ R(x)
  \labl tw  \end{equation}
for all $x\iN\L$. Second, when $(R{\circ}\omega,V)$ happens to be 
isomorphic to $(R,V)$ and when in addition
$(R,V)$ is a weight module with \findim\ weight spaces -- for which the notion
of a (formal) {\em character\/} $\chii_V\eq{\rm tr}_V q^H$ 
is available -- then one also has a `twisted' or `twining' analogue 
  \be  \chii^\omega_V:= {\rm tr}_V^{}\, \Tau_\omega q^H  \labl tc  \end{equation}
of $\chii_V$. The {\em twining characters\/} \erf tc are generalized 
character-valued indices. Their coefficients in an expansion in powers of the 
formal variables $q$ are elements in an extension of $\zet$ by the eigenvalues 
of $\omega$, i.e.\ in a cyclotomic field when $\omega$ has finite order.

\sect{Enhancement of symmetries}

Symmetry enhancement in \qft\ refers to situations where one starts from a known 
model with symmetry $\chir$ and attempts to obtain a new  model with extended
symmetry $\chir\Xxt\,{\supset}\,\chir$. A general mechanism for doing so consists
in promoting relatively local (charged) fields to observables. 
In CCFT, this amounts to include suitable intertwining operators into the \voa\
\cite{dolm,lI4'}. As already announced, the extensions that we consider here 
are extensions by simple current fields. That is, there is some simple current 
group $\calg$ of the original CFT model such that
  \be  \calh\XXt_{\vac\xxt} = \mbox{\Large$\bigoplus$}^{}_{\J\in\calg} \calh_\J
  \end{equation}
is the underlying vector space of the chiral \alg\ of the extended model.
The chiral \Cft\ model with enhanced 
symmetry is called the {\em simple current extension\/} of the original one.

A complete understanding of the extended chiral \alg\ $\chir\Xxt$ and its \rep\ 
theory is in general not yet available. Among the ingredients of the extended 
model that one would like to establish are in particular:\,%
 \footnote{~Those notions appearing here that were not yet introduced will 
 be explained later on.}
\Nxt 
The label set $\sus\Xxt\eq\{\mu\xxt\}$ of $\chir\Xxt$-sectors.
\nxt  
The decomposition $\calh\XXt_{\mu\xxt}\eq\bigoplus_i\!\calh_{\mu_i}$
of $\chir\Xxt$-sectors viewed as \chir-sectors.
\nxt
The fusion rules ${\rm N}\Xxt$ of the $\chir\Xxt$-model.
\nxt 
The characters $\chii\XXt_{\mu\xxt}$ of the $\chir\Xxt$-model and their modular 
S-transformation matrix.
\nxt  
The \cb s and \corfu s of the $\chir\Xxt$-model.

\vskip.3em
To become more concrete, let us first collect a few basic facts 
\cite{scya,scya6} that are needed to understand simple current extensions. 
We start with a fusion ring containing a non-trivial group of simple currents. 
But the perspective is somewhat different from the exposition in section 2, in
that the Verlinde conjecture (see below) is built in from the beginning; thus 
the matrix $S$ is provided by the modular transformation of characters,
and then the fusion rules are {\em defined\/} through $S$ by formula \erf1v. 

The input information needed for the extension then consists of the following:
\Nxt
A set $\{ {\chii_\mu} \}$ \,($\mu\iN I$, $|I|\,{<}\,\infty$)
of functions of a complex variable $\tau$, convergent in the upper half-plane and
forming a basis of a unitary module $V$ over \slz\ for which the generator
$S$ (implementing $\tau\,{\mapsto}\,{-}1/\tau$) is symmetric and the generator 
$T$ (implementing $\tau\,{\mapsto}\,\tau{+}1$) is diagonal.
\nxt
A vacuum label $\vac\iN I$, satisfying $S_{\vac,\mu}\iN\reals_{>0}$
for all $\mu\iN I$. 
\nxt
An involution $C{:}\;\mu\,{\mapsto}\mu^+$ on $I$ leaving $\vac$ fixed and
satisfying $CS\eq S^*$ and $CT\eq TC$.
\nxt
A subset $\calg\,{\subseteq}\,I\,$ such that
$S_{\J,\vac}\eq S_{\vac,\vac}$ and $T_\J\eq T_\vac$ for all $\J\iN\calg$
\nxt
Numbers $\N\lambda\mu\nu$, defined by \erf1v through the modular
transformation matrix $S$.
\nxt
A product `$\star$' on $\complex^{|I|}$ defined by formula \erf1N.
\vskip.3em

Having these data at one's disposal, many of the goals in the 
list given above can be reached. Concretely, one first proves \cite{scya6}:
\Nxt
$\complex^{|I|}$ endowed with the product `$\star$' is a fusion \alg\
$\cal A$.
\nxt
$\calg$ is a finite abelian group w.r.t.\ `$\star$'
-- the group of units in the basis $I$ of $\cal A$.
\nxt
$\calg$ organizes $I$ into orbits $[\mu]\,{:=}\,\{\J\mu\,|\,\J{\iN}\calg \}$
with $\J\mu\,{:=}\,\J\,{\star}\,\mu\iN I$. The length of the orbits is
$\ell_\mu\eq|\calg|/|\cals_\mu|$,
with $\cals_\lambda\,{:=}\,\{\J\iN\calg\,|\,\J\lambda\eq\lambda\}$ the
stabilizer of $\lambda$ in $\calg$.
\nxt
The combination
  \be  Z^{(\calg)}(\tau) := \sum_{\scs\mU: \quad \mu\in I,\ \ \atop
  \scs T_{\J\mu}=T_\mu\ \forall\J\in\calg}\!\!\!
  \llb\, |\cals_\mu|\cdot \mbox{\LARGE$|$}\!\! \sum_{\J\in\calg/\cals_\mu}\!
  \chii_{\J\mu}(\tau) \mbox{\LARGE$|$}^2 \,\lrb  \labl21
  \end{equation}
is invariant under the \pslz\,-action $\tau\,{\mapsto}\,\frac{a\tau+b}{c\tau+d}$.
\vskip.3em

Furthermore, the following interpretation has been shown to be self-consistent
in \cite{fusS6} and, when combined with uniqueness results \cite{brug2,muge6}
on the modularisation of premodular categories, can be established rigorously:
\Nxt
$Z^{(\calg)}$ as given by \erf21 is the `diagonal modular invariant' 
  \be  Z^{(\calg)}(\tau) = Z\Xxt(\tau) := \mbox{\Large$\sum$}^{}
  _{\mu\xxt\in\sus\xxt} |\chii_{\mu\xxt}(\tau)|^2_{\phantom I}  \end{equation}
for the simple current extension.
\nxt 
The label set $\sus\Xxt$ consists of all equivalence classes of pairs $\Mu$ with
$\mu\iN\sus$ obeying $T_{\J\mu}\eq T_\mu$ and $\psu$ a character of the 
{\em central stabilizer\/}
  \be  \calu_\mu := \{\J\iN\cals_\mu \,|\, F_\mu(\J{,}\J'){=}1
  \,\forall\, {\J'}{\in}\cals_\mu\}\;\subseteq\;\cals_\mu \,.  \labl33
  \end{equation}
Here $F_\mu$ is a certain alternating bihomomorphism on $\cals_\mu$ (see 
formula \erf39 below), and the classes $\Mu$ are with respect to the action 
$(\mu,\psu)\,{\mapsto}\, (\J\mu,\overline{F_\mu(\J,{\cdot}\,)\psu})$ for $\J\iN\calg$.
\nxt 
The functions $\chii\Xxt$ are given by
  \be  \chii\XXt_{[\mu,\psu]} = d_\mu{\cdot}\mbox{\Large$\sum$}^{}_
  {\J\in\calg/\cals_\mu} \chii_{\J\mu} \labl ce \end{equation}
with $d_\mu$ the square root of the embedding index of $\,\calu_\mu\,{\subseteq}
\,\cals_\mu$. These functions can be interpreted as the characters of simple 
modules $\calh\XXt_{[\mu,\psu]}$ over the extended VOA.\,%
 \footnote{~As an illustration of the appearance of the prefactors $d_\mu$, take
 the WZW model based on $\gb\eq\so(n_1)\oplu\so(n_2)\oplu\so(n_3)$ at 
 level 1 with $n_1,n_2,n_3\iN2\zetpluso{+}1$, extended by the simple current group
 $\{ (\vac{,}\vac{,}\vac),(\vac{,}\vv{,}\vv),(\vv{,}\vac{,}\vv),
 (\vv{,}\vv{,}\vac)\}\cong\zet_2\,{\times}\,\zet_2$, with $\vv$ referring to the
 vector \rep\ of $\so(n)$. This extension can also be understood as a conformal
 embedding into $\so(n_1{+}n_2{+}n_3)$ at level 1. There is a single fixed point
 $(\ss{,}\ss{,}\ss)$, where $\ss$ refers to the $\so(n)$ spinor \rep. The 
 analysis shows that $d\eq2$ for this fixed point. This is in agreement with the
 fact that the dimension $2^{(n_1+n_2+n_3-1)/2}$ of the irreducible spinor \rep\
 of $\so(n_1{+}n_2{+}n_3)$ is twice as large as the dimension
 $2^{(n_1-1)/2} 2^{(n_2-1)/2} 2^{(n_3-1)/2}$ of $(\ss{,}\ss{,}\ss)$.}
\nxt 
Thus the summands in the expression \erf21 are to be
read as $|\calu_\mu|\,{\cdot}\,|\chii\XXt_{[\mu,\cdot]}|^2$, i.e.\ to each orbit
$[\mu]$ there are associated $|\calu_\mu|$ many extended irreducible characters
$\chii\XXt_{[\mu,\psu]}$. Correspondingly, considered as an \chir-module,
the $\chir\Xxt$-module $\calh\XXt_{\mu\xxt}$ decomposes as
  \be  \calh\XXt_{[\mu,\psu]} = \mbox{\Large$\bigoplus$}^{}
  _{\J\in\calg/\cals_\mu} \complex^{d_\mu} {\otimes} \calh_{\J\mu} \,. \end{equation}
\nxT 
The functions $\chii\XXt_{[\mu,\psu]}$ \erf ce span a unitary \slz-module.
Their S-transformation matrix $S\Xxt$ is obtained by sandwiching certain 
matrices $S^\J$ between characters of central stabilizers:
  \be  S\XXt_{[\lambda,\psu_\lambda],[\mu,\psu_\mu]}
  = \,\Frac{|\calg|^{\phantom|}}{[\,|\cals_\lambda|\,|\calu_\lambda|\,
  |\cals_\mu|\,|\calu_\mu|\,{]}^{1/2}_{\phantom I}}
  \sum_{\J\in\calu_\lambda\cap\calu_\mu}
  \psu_\lambda(\J)\, \SJ_{\lambda,\mu}\, {\psu_\mu}(\J)^* \,.
  \labl so  \end{equation}
The matrices $\SJ$ occurring here are those that describe the modular 
S-transformation of the one-point chiral blocks $B_\J^{\sss(1)}$ with insertion
$\J$ on an elliptic curve. (The rows and columns of $\SJ$ are labelled by
$\{\mu\iN\sus\,|\,\cals_\mu\,{\ni}\,\J\}$, and $\SJ$ obeys the \slz\
relations as well as $(S^\J)^{\rm t}_{}\eq S^{\J^{-1}}\!$, $S^\vac\eq S$.)
\nxt 
$S\Xxt$ as given by \erf so is proven to be unitary and symmetric and to satisfy 
the \slz\ relations. With the help of the computer program {\tt kac} \cite{bs} 
it was also checked in a huge number of examples that \erf so produces
non-negative integers when inserted into formula \erf1v \cite{fusS6}.
\vskip.3em

Several further comments are in order: 
\Nxt 
There is evidence that
for WZW models (and similarly for coset models), up to a calculable phase 
$\SJ$ is the Kac\hy Peterson matrix of the {\em orbit \lie\/} that is
associated \cite{fusS3,furs}\,%
 \footnote{~This is based on the identity between twining characters for \g\
 and ordinary characters for the orbit \lie, which was established in 
 \cite{fusS3,furs} for Verma modules and simple highest weight modules.
 Results on other classes of modules can be found in \cite{nait11,kakk,kaNa},
 and related work in \cite{schW3,bofm,wend2}.}
to \g\ and $\J$ via a folding of the Dynkin diagram. Hence for these models
all $\SJ$ are known explicitly.
\nxt 
By general results on the cohomology of finite abelian groups,
the alternating two-cocycle $F_\mu$ is the commutator cocycle for some 
${\cal F}_{\!\mu}\iN H^2(\cals_\mu{,}{\rm U}(1))$, i.e.\ satisfies
$F_\mu(\J{,}\J')\eq{\cal F}_\mu(\J{,}\J')/{\cal F}_\mu (\J'{,}\J)$
\cite{bant6,muge6,fuSc112}. Thus the group \alg\ $\complex\,\calu_\mu$
is isomorphic to the center of the twisted group \alg\
$\complex_{{\cal F}_\mu}\!\cals_\mu$, implying e.g.\ that the inclusion
$\calu_\mu\,{\subseteq}\,\cals_\mu$ is of square index $d_\mu^2$, with
$d_\mu$ the dimension of the irreducible 
$\complex_{{\cal F}_\mu}\!\cals_\mu$-\rep s. (The occurrence of 
non-trivial cohomology classes ${\cal F}_\mu$ can be regarded
as a manifestation of the fact that in quantum theory symmetries 
are in general only realized projectively.)
\nxt 
It follows directly from the definition of $\SJ$ as the modular 
S-transformation matrix for the chiral blocks $B_\J^{\sss(1)}$ that
  \be  \Frac{S^{\K}_{\J\lambda,\mu}}{S_{\vac,\vac}}
  = \Frac{S_{\mu,\J}}{S_{\mu,\vac}}
  \cdot \F{\J\lambda}\lambdap\lambdap \K{\J\lambda}\Jp
  \cdot \Frac{S^{\K}_{\lambda,\mu}}{S_{\vac,\vac}} \,,  \labl1f  \end{equation}
where $\F{\J\lambda}\lambdap\lambdap \K{\J\lambda}\Jp$
is an entry of a so-called fusing matrix {\sf F}, i.e.\ a $6j$-symbol of the 
modular tensor category of the original \Cft.
One can also show \cite{fusS6} that
  \be  \SJ_{\J'\lambda,\mu} = (T_\mu/T_{\J'\mu})\,F_\mu(\J,\J')
  \,\SJ_{\lambda,\mu} \,,  \end{equation}
which by comparison with \erf1f tells us that
  \be  F_\lambda(\J,\K) = \F{\J\lambda}\lambdap\lambdap\K{\J\lambda}\Jp \,.
  \labl39 \end{equation}
\nxT 
In connection with \bc s also simple currents of conformal weight $\Delta
\iN\zet{+}\frac12$ turn out to be relevant. In that case $F_\mu$ 
is no longer a commutator cocycle. But after combining it with so-called 
discrete torsion \cite{krSc}, one arrives again at a cohomological 
interpretation \cite{fhssw}.  
\nxt 
Recall the relation \erf gh between the simple currents of the \wzwm\ 
based on \g\ and the center of the Lie group \G\ (for which 
${\rm Lie}(\G)$ is the compact real form of the horizontal sub\alg\ 
\gb\ of \g). It says that (except for $\gb\eq E_8$ and level 2) 
  \be  \calg(\g) \cong {\cal Z}(\G) \cong \lwv/\Lv \labl ll 
  \end{equation}
with $\lwv$ and $\Lv\,{\subseteq}\,\lwv$ the coweight and coroot lattices of 
\gb, \resp. (${\cal Z}(\G)$, and hence $\calg(\g)$, is also naturally
isomorphic to the maximal normal abelian subgroup of the group of symmetries
of the Dynkin diagram of the affine \lie\ \g.) It is therefore not surprising 
that the matrix $S\Xxt$ appears in the Verlinde formula for non-simply connected 
quotients $\tilde\G$ of $\G$. For some of the matrix elements $S\Xxt_{\vac\xxt
,\mu\xxt}$ the prediction \erf so was checked in \cite{beau3} for 
$\tilde\G\eq{\rm PGL}(n{,}\complex\,)$ with $n$ prime.
\nxt 
Among the labels $\mu\xxt\eq\Mu$ associated to a given orbit $[\mu]$ there 
is no distinguished one. Accordingly, the degeneracy labels $\psu$ should 
better not be regarded as elements of the character
group $\calu_\mu^*$, but rather as elements of the torsor over that group. 
Insisting on the description in terms of $\calu_\mu^*$, certain non-canonical 
basis choices are implied. A formulation free of such non-canonical 
choices can be achieved \cite{bai} by expressing $\mu\xxt$ through 
suitable idempotents in morphism spaces, similarly as in \cite{brug2,muge6}.

\sect{Simple current automorphisms of block \alg s and chiral blocks}

The modular tensor category summarizes the basic quantities of a chiral CFT
model, such as the set $\sus$ of sector labels, the fractional part of 
conformal weights $\Delta$, the fusion rules with their diagonalizing matrix 
$S$, and also the fusing matrices {\sf F}. It contains 
more sophisticated information as well,
in particular about the {\em \cb s\/} of the model on arbitrary curves $\tC$. 
The chiral (or {\em conformal\/}) blocks $B_\muv(\tC_g)$ for a genus-$g$ 
curve $\tC_g$ with distinct marked points $p_1,p_2,...\,,p_m$ and 
associated sector labels $\mu_1,\mu_2,...\,,\mu_m\iN\sus$ are specific 
linear forms $B_\muv{:}\ \calhmv\,{\to}\,\complex\,$, where
(tensor products are over $\complex$\,)
  \be  \calhmv := \calh_{\mu_1}\,{\otimes}\,\calh_{\mu_2}\,{\otimes}\,\cdots
  \,{\otimes}\,\calh_{\mu_m} \,. \end{equation}
Namely, they form the space of co-invariants \wrtt action of the {\em block 
\alg\/} $\apcg$ that expresses the chiral symmetries globally on $\tC_g$:
  \be  B_\muv(\tC_g) = [\calhmv]_\apcg^{\phantom I} \,.  \end{equation}
Alternatively, by duality one may think of the \cb s as the invariants 
$\llb(\calhmv)^*_{\phantom I}\lrb^\apcg_{\phantom o}$
in the algebraic dual of $\calhmv$. In physical terms, the invariance
condition amounts to imposing the {\em Ward identities\/} for the symmetries.

For WZW models, the block \alg\ is a \lie, which is the tensor product 
\cite{tsuy,Ueno,beau}
  \be  \apc = \gb \otimes \fpc  \labl gf  \end{equation}
of the horizontal sub\alg\ \gb\ of the affine \lie\ \g\ with the associative
\alg\ $\fpc$ of meromorphic functions on $\tC$ whose singularities are at
most finite order poles at the marked points $p_s$. The action of the
\alg\ \erf gf on $\calhmv$ follows via Laurent expansion of the functions 
in $\fpc$ in local coordinates around the $p_s$.  

The corresponding prescription in the general case is considerably more 
involved \cite{freN5}.
One introduces a bundle ${\cal A}$ over $\tC$ whose fibers are isomorphic 
to the vector space $\calh_\vac$ of the VOA, and requires that for every 
$m$-tuple of vectors $v_i\iN\calh_{\mu_i}$ and every $w\iN\calh_\vac$ the 
section $B_\muv(v_1\ot\cdots\ot\,{\cal Y}(w,{\cdot}\,)v_i\,\ot\cdots\ot v_m)$
(with ${\cal Y}$ a section of ${\cal A}^*$
that takes over the role of the state-field correspondence) in the restriction 
of ${\cal A}^*$ over local disks around the points $p_s$ can be extended 
to a global holomorphic section on the punctured curve $\tC_g{\setminus}\vec p$.
In view of this general prescription, it is actually a non-trivial statement
that in \wzwm s it is sufficient to impose only the Ward identities coming 
from the affine \lie\ \g.

\smallskip

The following properties of the block spaces $B_\muv(\tC_g)$
have been proven in the WZW case \cite{tsuy} and are expected to hold for
arbitrary (rational, and partly even for non-rational) CFT models:
\nxt
For fixed $\vec p\eq(p_1,...\,,p_m)$ and fixed moduli of $\tC_g$,
$B_\muv(\tC_g)$ is a \findim\ vector space. 
\nxt
For fixed $\muv\eq(\mu_1,...\,,\mu_m)$ and fixed genus $g$ these spaces 
fit together to the total space of a vector 
bundle $\calbg_\muv$ of finite rank over the moduli space $\calm_{g,m}$
of complex curves of genus $g$ with $m$ ordered marked points. In particular, 
the dimension of $B_\muv(\tC_g)$ depends only on $g$ and on $\muv$, but neither 
on the moduli of the curve $\tC_g$ nor on the positions of the marked points.
\nxt
The dimensions of the genus zero 3-point \cb s 
$B_{\lambda,\mu,\nuP}^{\sss(0)}$ are given by the fusion rules
$\N\lambda\mu\nu$ \erf1N, and those for higher genera and\,/\,or more 
insertions can be expressed through the fusion rules by simple factorization
prescriptions. Together with \erf1v, this results in the formula
  \be  {\rm rank}\, \calbg_\muv = \sum_{\nu\in\sus} |S_{\vac,\nu}|^{2-2g}_
  {\phantom i}\, \prod_{s=1}^m \, \frac{S_{\mu_s,\nu}}{S_{\vac,\nu}}
  \,,   \labl vv \end{equation}
where $S$ is the matrix that according to \erf1v diagonalizes the \furu s.
\nxt
Each bundle $\calbg_\muv$ comes equipped with a pro\-jectively flat connection,
the \kzc. This implies a projective action of the mapping class group on
$B_\muv(\tC_g)$. In the case of an elliptic curve and with
a single `insertion' $\mu_1\eq\vac$, this furnishes a unitary
projective representation of the modular group \pslz\ on 
the 1-point chiral blocks $B_\vac^{\sss(1)}\,{\cong}\,\complex^{|\sus|}_{}$.
\nxt
In a natural basis, the generator 
S${:}\;\tau\,{\mapsto}\,{-}1/\tau$ of \pslz\ is represented by a symmetric 
matrix $S$.
This matrix coincides with the diagonalizing matrix that appears
in \erf vv, justifying that we already used the same symbol for both of them. 
\nxt
The {\em Verlinde conjecture\/} states that this matrix $S$ also coincides
with the matrix describing the modular S-transformation $\tau\,{\mapsto}\,{-}1/\tau$
on the characters $\chii_\mu$.
This conjecture is proven (see e.g.\ \cite{beau,falt,tele,fink,sorg}) for 
WZW models for which the modular S-matrix is given by the Kac\hy Peterson formula. 

\smallskip

There is no reason to expect that the vector bundles $\calbg_\muv$
are irreducible. Indeed any \auto\ of $\calbg_\muv$ constitutes trivially a
source for reducibility, namely into a direct sum of its eigenspaces.
Such \auto s can e.g.\ be inherited from suitable \auto s of the relevant
block \alg\ $\apc$. What is non-trivial is the observation that such \auto s 
of $\apc$ arise naturally, namely as a consequence of the presence of 
{\em simple currents\/}.

This can be made explicit for \wzwm s, for which $\apc$ takes the form \erf gf.
Then there are \auto s of $\apc$ coming from simple current \auto s of the 
underlying affine \lie\ \g. They are constructed as follows \cite{fuSc8}. For 
$g\eq0$ and any $m\,{\ge}\,2$ one has a group of \auto s of \g\ labelled by the set
  \be  \Gamma_{\!\!{\rm w}} := \{ (\nub_1,\nub_2,...\,,\nub_m) 
  \iN(\lwv)^m \,|\, \mbox{$\sum_{s=1}^m$}\nub_s\eq0\}  \labl gw \end{equation}
($\lwv$ the coweight lattice of \gb). These \auto s $\ommuv$ depend in addition 
on a sequence of pairwise distinct complex parameters $z_s$ ($s\eq\onetom$), one 
of which, say $z_1$, is singled out. $\ommuv$ acts on the canonical central element 
$K$ and on the elements $H^i\ot f$ and $\eal\ot f$ of $\g\eq\complex(\!(t)\!)
{\oplus}\complex K$\, ($f\iN\complex(\!(t)\!)$ 
and $\{H^i\,|\,i\eq1{,}2{,}...,{\rm rank}(\gb)\}{\cup}$\linebreak[0]$\{\eal\,|\,
\bar\alpha\;{\rm a}\;\gb\mbox{-root}\}$ a Cartan\hy Weyl basis of \gb) as
  \be  \bearl  \ommuv(H^i\ot f):= H^i\ot f
  + K \, \sum_{s=1}^m \nub_s^i\, \Res(\fiss\,f) \,, \\{}\\[-.4em]
  \ommuv(\ebe\ot f):= \ebe\ot f \,{\cdot} \prod_{s=1}^m
  (\fiss)_{}^{-(\nub_s,\betab)} \,, \qquad  \ommuv(K) := K \eear \labl Au
  \end{equation}
with
  \be  \fiss(t) := \llb t+(z_1\,{-}\, z_s)\lrb_{}^{-1} \,.  \labl fs
  \end{equation}
Here $\Res$ denotes the residue of a formal Laurent series in $t$, and
we employ the notation $f^\ell$ for the function with values
$f^\ell(t)\eq(f(t))^\ell$. The {\em mul\-ti-shift 
automorphisms\/} \erf Au of \g\ are close relatives of the
ordinary shift automorphisms -- studied e.g.\ in \cite{frha,goom,levw} --
that can be recovered from formula \erf Au by setting $\nub_s\eq0$ for $s\nE 1$
(which does not belong to $\Gamma_{\!\!{\rm w}}$, though).

The set $\Gamma_{\!\!{\rm w}}$ \erf gw is an abelian group \wrt component-wise 
addition, isomorphic to $(\lwv)^{m-1}$. The \auto s \erf Au form a group 
isomorphic to $\Gamma_{\!\!{\rm w}}$.
Further, the function $\fiss$ \erf fs  can be recognized as the local expansion
at $p_1$ of the function $\varphi_{(s)}\iN{\cal F}_{\vec p,\Pe}$ defined by
  \be  \varphi_{(s)}(z) := (z-z_s)_{}^{-1} \,,  \end{equation}
where $z$ is a quasi-global coordinate on $\Pe$ and $z_s\eq z(p_s)$. It follows
\cite{fuSc8} that to every $\vec\nub\iN\Gamma_{\!\!{\rm w}}$ there is also
associated an \auto\ $\Ommuv$ of the block \alg\ $\gb{\otimes}\fpe$, acting as
  \be  \Ommuv(H^i\ot f) = H^i\ot f \,, \qquad \Ommuv(\ebe\ot f)
  = \ebe\ot\, f \cdot \prod_{s=1}^m (\varphi_{(s)})_{}^{-(\nub_s,\betab)}
  \,.  \labl aU  \end{equation}
Up to central terms the local expansions of the \auto\ \erf aU at the marked points 
reproduce the \auto s \erf Au of \g.\,%
 \footnote{~The central terms are needed in order to have an \auto\
 of the affine \lie\ rather than only of the corresponding loop \alg. Accordingly,
 the centers of the $m$ copies of \g\ associated to the marked points $p_s$
 must be identified; then upon summing over all insertion points these terms
 cancel by the residue theorem.}

Every \auto\ $\ommuv$ of \g\ induces a permutation $\sigma_{\vec{\bar\nu}}^*$ of 
the label set $I$. Further, for every $\vec\nub\iN\Gamma_{\!\!{\rm w}}$ and 
every $\muv\iN I^m_{}$ there is an induced map
  \be  \Tau_{\vec\nub}:\quad \calhmv\,{\to}\,\vec{\cal H}_{\Ommuvs\muv} \,,
  \labl im \end{equation}
unique up to a scalar, which is a twisted intertwiner in the sense of \erf tw. 
This descends to a map $\Tau^*_{\vec\nub}$ from the \cb s $B_\muv$ to 
$(\vec{\cal H}_{\Ommuvs\muv})^*$,
which is in fact an isomorphism to $B_{\Ommuvs\muv}$.  
The map on the co-invariants obtained this way from any {\em inner\/} \auto\
of the block \alg\ is just the identity. Thus for the action on \cb s
we are interested in the outer \auto\ class\ of $\Ommuv$ \erf aU.
Now $\Ommuv$ is outer iff at least one of the vectors $\nub_s$ in the 
tuple $\vec\nub\iN\Gamma_{\!\!{\rm w}}$ is not a coroot, so the group of 
outer \auto\ classes is the factor group
  \be  \Gamma_{\!\!{\rm out}} := \Gamma_{\!\!{\rm w}} / \Gamma 
  \quad\ {\rm with}\quad\ 
  \Gamma := \{ (\bub_1,\bub_2,...\,,\bub_m)
  \iN(\Lv)^m \,|\, \mbox{$\sum_{s=1}^m$}\bub_s\eq0\} \,.  \labl gg
  \end{equation}
$\Gamma_{\!\!{\rm out}}$ is a finite abelian group, obeying (compare formula \erf ll)
  \be  \Gamma_{\!\!{\rm out}} \cong (\lwv/\Lv)^{m-1}_{}
  \cong (\calg(\g))^{m-1}_{} \,,  \end{equation}
with $\calg$ the group of simple currents of the \wzwm. Thus we can view
the elements of $\Gamma_{\!\!{\rm out}}$ as $m$-tuples of simple currents 
$\J_s$ subject to $\J_1{\star}\J_2{\star}\cdots{\star}\J_m\eq\vac$,
and accordingly write the map $\Tau^*_{\vec\nub}$ as 
  \be  \Tau^*_{\J_1,...,\J_m}:\quad B_\muv \to B_{\Ommuvs\muv} \,.
  \labl mi \end{equation}

When $\J_s\iN\cals_{\mu_s}$ for all $s\eq\onetom$, then the map
$\Tau^*_{\J_1,\ldots, \J_m}$ is an {\em auto\/}morphism of $B_\muv$, and 
in particular one can study its trace. So
for every $\muv\iN I^m$ the group $\cals_\muv\,{:=}\,\Gamma_{\!\!{\rm out}}
{\cap}\llb\mbox{\LARGE$\times$}\!{}_{i=1}^{m_{\phantom i}} \cals_{\mu_s}\lrb$
acts on $B_\muv$ by the maps $\Tau^*_{\J_1,...,\J_m}$. Since the twisted 
intertwining property determines the maps $\Tau_{\vec\nub}$ \erf im, and 
hence also the maps $\Tau^*_{\J_1,...,\J_m}$ \erf mi, 
only up to a phase, our construction furnishes in general only a projective 
\rep\ of $\cals_\muv$ on $B_\muv$. We denote by $\calf_\muv$ the class in 
$H^2(\cals_\muv,{\rm U}(1))$ that specifies the projectivity. 

The cohomology class $\calf_\muv$ depends on $\muv$, as well as a priori on 
the insertion points 
$\vec p$. But one can choose \cite{fuSc8} these phases in such a manner 
that $\calf_\muv$ does in fact not depend on $\tau$ and $\vec p$. Also, with 
this choice the map $\Tau^*_{\J_1,...,\J_m}$ has finite order, and the
choice is compatible with the \kzc\ on $\calb_\muv$.
Thus for every element $\vec\J\iN\cals_\muv$ we have a finite order \auto\ of 
the chiral block {\em bundle\/} $\calb_\muv$ over the moduli space
$\calm_{0,m}$. We will use the same symbol $\Tau^*_{\J_1,...,\J_m}$ for
this bundle \auto\ as for the vector space \auto s of the individual
fibers $B_\muv$ of the bundle.

We are interested in the ranks of the sub-bundles that are invariant
under the action of $\Tau^*_{\J_1,...,\J_m}$; they depend on the tuple 
$\vec\J$ of simple currents. A hint on how this dependence looks like can 
be obtained by inspecting the Verlinde formula for simple current 
extensions of WZW models, in which the matrix \erf so appears.
This leads to a concrete conjecture
for the ranks of the sub-bundles in the case of genus 0, which can be
extended to higher genera by analogy with the Verlinde formula. 

We first formulate a conjecture for the cohomology class of $\calf_\muv$, 
in terms of its commutator cocycle $F_\muv$: We propose that
$F_\muv$ is the product of the commutator cocycles of the factors:
  \be  F_\muv = \mbox{\Large$\prod$}_{i=1}^m F_{\mu_i} \,,  \end{equation}
   with $F_{\mu_i}$ given by \erf39.
It is then natural to introduce the `central stabilizer' (compare formula \erf33)
$\calu_\muv\,{:=}\,\{\vec\J\iN\cals_\muv\,|\, F_\muv(\vec\J{,}\vec\J'){=}1
\,\forall\, {\vec\J'}{\in}\cals_\muv\}$, for which our conjecture amounts to
  \be  \calu_\muv = \Gamma_{\!\!{\rm out}}
  {\cap}\llb\mbox{\LARGE$\times$}\!{}_{i=1}^{m_{\phantom i}} \calu_{\mu_s}\lrb
  \,.  \end{equation}
For each character $\psi$ of $\calu_\muv$, there is now an invariant sub-bundle
$B_\psi$ of the bundle of chiral blocks. (A word of warning is in order: The
correspondence of characters $\psi$ and sub-bundles is {\em not\/} canonical.\,%
 \footnote{~The requirement that $F_\muv$ be constant on $\calm_{0,m}$
 does not yet determine the phase choices for the maps
 $\Tau^*_{\J_1,...,\J_m}$ uniquely. Different allowed choices result in 
 different values for the traces. But one can show that this
 difference simply amounts to a relabelling of the eigenspaces.
 So again one deals with the torsor over the character group $\calu_\muv^*$
 rather than with $\calu_\muv^*$ itself.}
Our (conjectural) formula for the rank of these sub-bundles,
however, depends on precisely the same non-canonical choices.)

Rather than describing the ranks of the sub-bundles directly, it proves
to be more convenient to perform a Fourier transform over the central 
stabilizer $\calu_\muv$ and to give instead the traces
of the maps $\Tau^*_{\J_1,...,\J_m}$. Our conjecture is then formulated 
as follows.
\begin{quote}
 {\bf Conjecture 1:} The trace on $B_{\vec\mu}^{\sss(g)}$ of the twisted
 intertwiner $\Tau^*_{\J_1,...,\J_m}$ associated with the $m$-tuple
 $(\J_1,\J_2,...\,,\J_m)\iN\calu_\muv$ is
  \be  {\rm tr}_{\!B_{\vec\mu}^{(g)}}^{} \Tau^*_{\J_1,...,\J_m}
  = \sum_{\nu\in\sus} |S_{\vac,\nu}|^{2-2g}_{}
  \prod_{s=1}^m \, \frac{S^{\J_s}_{\mu_s,\nu}}{S_{\vac,\nu}} \,
  \,.  \labl44 \end{equation}
\end{quote}
Here $S^\J$ is the matrix introduced in \erf so, i.e.\ the modular
S-transformation matrix for the one-point chiral blocks $B_\J^{\sss(1)}$.
Also recall that in WZW models $S^\J$ is conjectured to coincide
(up to a calculable phase) with the Kac\hy Peterson matrix for the orbit
\lie\  that is associated to \g\ and $\J$; inserting this relation into
\erf44, one arrives at a numerical prediction for
${\rm tr}_{\!B_{\vec\mu}^{(g)}} \Tau^*_{\J_1,...,\J_m}$
that is as explicit as the Verlinde formula.

Note that the basic ingredients in the formula \erf44 are closely related to 
simple current extensions. Since that construction works for arbitrary CFT 
models, not only in the WZW case, the conjecture \erf44 can be extended to 
every (rational) CFT model, too. But already
in the WZW case one is still far from being able to prove the formula. 
In particular, since $S^\vac$ is the ordinary matrix $S$,
for $\J_1\eq\cdots\eq\J_m\eq\vac$, \erf44 reduces to the Verlinde formula 
\erf vv for the trace of $\Tau^*\,{=}\,{\rm id}$. So any proof of \erf44 
is a fortiori also a proof of the Verlinde formula itself. On the other hand, 
as is immediately checked, \erf44 is compatible with factorization.
Moreover, in the WZW case there is enormous numerical evidence. 
Namely, for very many cases it has been checked on the computer that the 
numbers obtained by Fourier transforming the traces \erf44 are non-negative 
integers, as required for the interpretation as dimensions of eigenspaces, 
even though they are obtained as complicated combinations of arbitrary\,%
 \footnote{~As a consequence of $|\sus|\,{<}\,\infty$ they are in fact
integers in a cyclotomic extension of $\mathbb Q$.} complex numbers.

Surprisingly, these numerical studies reveal that already the traces themselves
are integral (they may be negative, though), even when the order of the 
\auto\ exceeds two. So far no explanation of this empirical observation
is available. But there is some reminiscence with the interpretation of the
dimensions of block spaces as the Euler number of a suitable (BGG-like) complex,
with non-negativity implied by an additional acyclicity property \cite{tele}. 
One may hope that a similar description can be found for the traces as well.

\sect{Restriction of observables}

We now switch to our second main topic. The
basic idea is just the reverse of what we have done in section 3. One starts 
from a QFT model with observables $\chir$ and a group $G$ of \auto s of $\chir$,
and desires to construct the $G$-{\em orbifold\/} of the theory, by which one 
means some new QFT model whose observables are given by the fixed point 
sub\alg\ \chirg\ of \chir\ \wrt $G$. In CFT, where the observables are 
realized as a VOA, it is natural to demand that the orbifold is again a CFT 
model. This requires in particular that the fixed point \alg\ \chirg\ contains
a Virasoro element, and in fact this must coincide with the Virasoro element 
$\vvir$ of \chir. Note that this condition implies that none of the elements 
of $G$ can be a (non-trivial) simple current \auto.

In the orbifold model, the \chir-\rep s $R$ and $R^\sigma\,{\equiv}\,R\,{\circ}
\,\sigma$ with $\sigma\iN G$ become indistinguishable, and hence describe 
the same sector (or collection of sectors) of the orbifold. On the 
other hand, when $R^\sigma$ is isomorphic
to $R$ already as an \chir-\rep, then $R$, even when it is \chir-irreducible,
becomes reducible in the orbifold. Several distinct orbifold sectors are then
obtained by splitting such an \chir-sector into eigenspaces of the twisted 
intertwiner $\Tau_\sigma$
that implements $\sigma$ on the sector. But in addition there are 
also further orbifold sectors that are {\em not\/} obtainable by 
decomposing \chir-sectors. They are called {\em twisted\/} sectors.

Again, the goal is to express all quantities of interest of the new theory,
i.e.\ the orbifold, `in terms of' the original model. Due to the presence
of twisted sectors this can be quite hard. Fortunately, it can be 
expected that the characters of twisted sectors are obtainable via the modular 
S-transformation of (differences of) the characters of untwisted sectors.
This relationship provides a powerful tool for concrete calculations, and
we will assume that it is indeed satisfied. Once this assumption is made,
a rigorous derivation of the results collected below is possible \cite{bifs}.

Every $\sigma\iN G$ induces a permutation $\sigma^*$ of the label set 
$\sus$ such that $\calh_\mu^\sigma \cong \calh_{\sigma^*\!\mu}$ as 
\chir-modules; $\sigma^*$ is an \auto\ of fusion rules.  
{}From now on we restrict our attention to orbifolds by a $\zet_2$-group. 
Then one must only distinguish between length-two $\sigma^*\!$-orbits of 
\chir-sectors and `fixed points' $\mu\eq\sigma^*\!\mu$. One finds:
\Nxt
Every length-two orbit $\{\mu,\sigma^*\!\mu\}$ gives rise to a single untwisted 
orbifold sector, which we denote by
$(\mu,0,0)$ and whose character is $\chiO_{(\mu,0,0)}(\tau)\eq\chii_\mu(\tau)$. 
\nxt
Every fixed point $\mu$ yields two distinct untwisted orbifold sectors 
$(\mu,\psi,0)$, $\psi\iN\{\pm1\}$. Their characters read
  \be  \chiO_{(\mu,\psi,0)}(\tau)
  = \Frac12\, {\rm tr}_{\calh_\mu} \llb (\bfe\,{+}\,\psi\,\Tau_\sigma)\,
  q^{L_0-c/24} \lrb = \Frac12\, \llb \chii_\mu(\tau)
  + \psi\, \eta_\mu^{-1}\, \chie_\mu(2\tau) \lrb \,.  \end{equation}
Here we write $\eta_\mu^{-1}\chie_\mu(2\tau)$ for the twining character
$\chii_\mu^\sigma(\tau)$ because for suitable choice of the phases
$\eta_\mu$ the functions $\chie_\mu(\tau)$ introduced this way possess an 
expansion with integral powers of $q\eq\eE^{2\pi\ii\tau}$. 
Also, $\TO_{(\lambda,\psi,0)}\eq T_\lambda^{}$.
\nxt
In addition, each fixed point also gives rise to two twisted orbifold sectors, 
which we label as $(\mud,\psi,1)$ with $\psi\iN\{\pm1\}$, with characters 
  \be  \chiO_{(\mud,\psi,1)}(\tau)
  = \Frac12\, \llb \chiz_\mud(\Frac\tau2) + \psi\,
  (\Tz_\mud)^{-1/2}_{}\,\chiz_\mud(\Frac{\tau+1}2)\lrb \,.  \end{equation}
The labels $\mud$ appearing here
are in one-to-one correspondence with the fixed points $\mu$. But generically
this correspondence is not canonical, so that in particular we cannot dispense 
of using two different kinds of labels for the functions $\chie$ and $\chiz$.
\nxt
The functions $\chiz_\mud(\tau)$ possess an expansion with integral powers 
of $q$, too. The T-transformation of the twisted orbifold characters therefore reads
  \be  \chiO_{(\mud,\psi,1)}(\tau{+}1)
  = \psi\, (\Tz_\mud)^{1/2}_{}\, \chiO_{(\mud,\psi,1)}(\tau) \,, \end{equation}
so that we have $\TO_{(\lambdad,\psi,1)}\eq\psi\,(\Tz_\lambdad)^{1/2}_{}$.
\nxt
The twisted and untwisted sectors transform into each other under an
S-transfor\-ma\-ti\-on. For the functions $\chie$ and $\chiz$, this amounts to
the relations
  \be  \chie_\lambda(-\Frac1\tau)
  = \mbox{\Large$\sum$}^{}_\mud\, \Se_{\lambda,\mud}\,\chiz_\mud(\tau)
  \qquad{\rm and}\qquad  \chiz_\lambdad(-\Frac1\tau) = \mbox{\Large$\sum$}^{}_\mu
  \,\Sz_{\lambdad,\mu}\, \chie_\mu(\tau)  \end{equation}
with a unitary matrix $\Se$ and $\Sz_{\lambdad,\mu}\eq\eta_\mu^{-2}\,\Se_{\mu,\lambdad}$.
\nxt
It follows that the modular S-transformation matrix of the orbifold is given by
  \be \bearll
  \SO_{(\lambda,0,0),(\mu,0,0)} = S_{\lambda,\mu} + S_{\lambda,\gs\!\mu} \,,\quad &
  \SO_{(\lambda,\psi,0),(\mu,\psi',0)} = \Frac12\, S_{\lambda,\mu} \,,
  \\{}\\[-.8em] 
  \SO_{(\lambda,0,0),(\mu,\psi,0)} = S_{\lambda,\mu}  \,, &
  \SO_{(\lambda,\psi,0),(\mud,\psi',1)} =
  \Frac12\, \psi\eta_\lambda^{-1} \Se_{\lambda,\mud} \,,
  \\{}\\[-.8em] 
  \SO_{(\lambda,0,0),(\mud,\psi,1)} = 0  \,,  &
  \SO_{(\lambdad,\psi,1),(\mud,\psi',1)} = \Frac12\,\psi\psi'\,P_{\lambdad,\mud}
  \end{array} \labl os \end{equation}
with
  \be  P := (\Tz)^{1/2}_{} (\Se)^{\rm t}_{} (\eta^{-1}\Te)^2_{} \Se (\Tz)^{1/2}_{}
  \,.  \end{equation}
\nxt
$\Sz\Se$ is an order-2 permutation which up to sign factors coincides with $P^2$,
while $\eta^{-1}\Se\Sz\eta$ is an order-2 permutation up to sign factors. 

\smallskip

To compute the functions $\chie_\mu$ and their modular transformations, 
sufficiently detailed \rep\ theoretic information is needed. Therefore we now 
specialize again to the WZW case. Then we are dealing with \auto s of 
an untwisted affine \lie\ \g, which must preserve the (Sugawara) Virasoro \alg.
Every such map comes from an \auto, for brevity to be again denoted by $\sigma$, 
of the horizontal sub\alg\ $\gb$ of \g. When $\sigma$ has finite order, then 
in a suitable Cartan\hy Weyl basis of \gb\ it can be written as
  \be  \sigma^{\phantom|} = \omego \circ \sigmas \,,  \end{equation}
with $\omego$ a diagram automorphism and $\sigmas\eq\exp(2\pi\ii\, {\rm ad}_
{H_\ss}^{})$ an inner automorphism. $H_\ss\,{\equiv}\,(\ss,H)$ is an element 
of the Cartan sub\alg\ satisfying $\sigma(H_\ss)\eq H_\ss$. (The 
automorphisms of \findim\ simple \lie s have been classified \cite{ONvi}. 
For a list of all order-2 automorphisms see e.g.\ table 1 of \cite{bifs}.)

Further, we know that $\sigmas^*\eq\id$, hence $\sigma^*\eq\Omego^*$, and the
phases $\eta_\mu$ are given by $\eta_\mu\eq \exp(2\pi\ii(\ss,\mu))$. Finally, 
the twining characters $\chii_\mu^\sigma(\tau)\eq\eta_\mu^{-1}\chie_\mu(2\tau)$
coincide with ordinary characters of the orbit \lie\ $\g\ol\sigma$ that is 
associated to \g\ and $\sigma$ and hence are known very explicitly. The orbit 
\lie s arising here turn out to be twisted affine \lie s. More precisely, 
$\g\ol\sigma$ isn't isomorphic to any untwisted \alg\ (i.e.\ is `genuinely 
twisted') iff $\sigma$ is outer. For inner $\sigma$, where $\g\ol\sigma$ is 
isomorphic to an untwisted \alg, there is a {\em canonical\/} correspondence 
between fixed points $\mu$ and twisted sector labels $\mud$.

Next, we exploit the behavior of the functions $\chie_\mu$ 
\wrtt Weyl group $W$ of \g, \resp\ \wrtt subgroup
  \be  W_{\!\omego} := \{ w\iN W\,|\,w\,{\circ}\,\sigma^*\eq\sigma^*{\circ}w \} \,,
  \end{equation}
which is isomorphic, as a Coxeter group, to the Weyl group of the orbit \lie\
$\g\ol\omego$ \cite{furs}. The functions $\chie$ and $\chiz$ can then be shown 
to be quotients of alternating $W_{\!\omego}$-sums of Theta functions. Except for
$A_{2n}$ with outer \auto, the result reads
  \be  \chie_\mu(2\tau) = \chii^{\omego}_\mu[0,\ss](\tau)  \qquad{\rm and}\qquad
  \chiz_\mud(\Frac\tau2) = \chi^{\omego}_\mud[\ss,0](\tau)  \labl01 \end{equation}
with $\mu\iN\Mw_{\!\!\!\!\!\omego}/(\bar W_{\!\omego}{\Ltimes}\,h\Lv_{\!\!\omego})$
and $\mud\iN\Lw_{\!\!\!\!\!\omego}/(\bar W_{\!\omego}{\Ltimes}\,h\MV_{\!\!\omego})$.
Here $\chii[\ss_1,\ss_2]$ are shifted characters defined by means of
shifted Theta functions, and the relevant lattices are $\Lv_{\!\!\omego}\eq
\{\sum_i n_i\alpha^{(i)\Vee}\,|\,n_{\dot\sigma i}\eq n_i\iN\zet \}$ and
$\MV_{\!\!\omego}\eq\{ \sum_i n_i\alpha^{(i)\Vee}\,|\,\ell_i n_i\iN\zet,\, 
n_{\dot\sigma i}{=}n_i \}$, with $\ell_i$ the orbit lengths \wrtt symmetry of
the Dynkin diagram of \gb\ that corresponds to $\omego$.

In the exceptional case $\gb\eq A_{2n}$ with $\sigma\eq\omego\eq\sigma_{\rm c}$
(`charge conjugation'), owing to the minus sign that (only) in this case appears in
the transformation $\omego(E^\theta)\eq{-}E^\theta$ of the generator associated to 
the highest \gb-root, the formulas \erf01 get replaced by
  \be  \chie_\mu(2\tau) = T_\mu^{-1/2}
  \chii^{\omego}_\mu[0,0](\tau{+}\mbox{$\frac12$})
  = \chii^{\omego}_\mu[0,\ss_\circ](\tau)\,, \quad  
  \chiz_\mud(\Frac\tau2) = \chi^{\omego}_\mud[\ss_\circ,0](\tau) \,, 
  \end{equation}
where $\ss_\circ$ is the $A_{2n}$-weight 
$\ss_\circ\eq\frac14\,(\Lambda_{(n)}{+}\,\Lambda_{(n+1)})$.

It is then not too difficult to read off the S- and T-transformations. One finds
  \be  \Se_{\lambda,\mud} = S^{\omego}_{\lambda,\mud}  \qquad{\rm and}\qquad
  \Sz_{\mud,\lambda} = \eta_\lambda^{-2}\, S^{\omego}_{\lambda,\mud}
  \labl51 \end{equation}
with $S^{\omego}_{\phantom i}$ the S-matrix of the orbit \lie\ $\g\ol\omego$,\,%
 \footnote{~This matrix expresses the characters of integrable modules over
 the twisted affine \lie\ $\g\ol\omego$ at $-1/\tau$ in terms of the characters
 at $\tau$ of another twisted affine \lie\ \cite{KAc3}.
 In some cases that \alg\ coincides with $\g\ol\omego$, but then on the
 S-transformed side one deals with a different set of representations. \label3} 
and
  \be  \Tz_\mud = \left\{\bearll
  \eE^{2\pi\ii k(\ss_\circ,\ss_\circ)} \eE^{2\pi\ii(\mu,\ss_\circ)}\,
  {(T_\mu)}^{1/2} & {\rm for}\ \gb\eq A_{2n},\; \omego\eq\sigma_{\rm c}\,,
  \\[.12em] \eE^{2\pi\ii k(\ss,\ss)}\, \eE^{2\pi\ii(\mud,2\ss)}\,
  {(T^{\omego}_\mud)}^2 & {\rm else} \,. \eear \right.  \end{equation}

\smallskip
  
To conclude this section, we remark that an action of the orbifold \auto s
$\sigma$ by twisted intertwiners $\Tau_\sigma$ can be defined on the 
\chir-modules $\calh_\mu$, and correspondingly there is an action on the 
spaces of chiral blocks of the \chir-theory. Inspection of the fusion rules of
the orbifold model \cite{bifs}, as obtained from the modular S-matrix 
$S^\calo$ \erf os by the Verlinde formula \erf1v, then leads to a conjecture
on the traces of the action of $\sigma^{\otimes m}$ on chiral blocks, 
similar to the formula \erf44: 
\begin{quote}
 {\bf Conjecture 2:} 
 For every \auto\ $\sigma$ of \L\ that preserves the Virasoro element
 $v_{\sss\rm Vir}$, the trace of the induced map on $B_{\muv;\Pe}$ is
  \be  
  {\rm tr}_{B_{\muv;\Pe}}^{} \llb \Tau_{\sigma,\sigma,...,\sigma}^* \lrb
  = \sum_{\dot\kappa} {|S_{\dot\kappa,\vac}^{\sss(0)}|}^2\,
  \prod_{i=1}^m \frac{S_{\dot\kappa,\mu_{\mbox{$\sss i$}}}^{\sss(0)}}
  {S_{\dot\kappa,\vac}^{\sss(0)}} \,,  \labl2C \end{equation}
 The matrix $S^{\sss(0)}$ appearing here 
 is given as in formula \erf51, i.e.\ coincides with the Kac\hy Pe\-ter\-son
 matrix of $\L\eq\g$ for $\sigma$ inner, and with the Kac\hy Peterson
 matrix of the orbit \lie\ $\g\ol\omego$ if $\sigma$ is outer.
 (It appears in the modular S-matrix of the orbifold, see \erf os.)
 \end{quote}
The conjecture is formulated for WZW models, but just like in the case of 
formula \erf44 it originates from structures present in every rational 
CFT model and hence can again be extended to arbitrary such models.
Note that in the present case the Fourier transformation between
the dimensions of eigenspaces and the traces \erf2C is to be performed 
\wrtt full cyclic group generated by $\sigma$. No projectivity needs to
be accounted for, since one and the same \auto\ $\sigma$ is used at each 
insertion point and the second cohomology of cyclic groups is trivial.

\sect{Full \cft}

We now turn our attention to what in contradistinction to chiral CFT
one calls {\em full\/} CFT. It is full CFT that is relevant to many 
applications, e.g.\ in statistical mechanics, condensed matter physics and 
string theory. (CCFT, on the other hand, appears in the description
of the (fractional) quantum Hall effect.)
The world sheet for full CFT is a real two-dimensional manifold $X$ with a 
conformal structure. $X$ may have a boundary, and it need not be orientable. 
Even when it is orientable, it does not come with a canonical choice of orientation.

In order to exploit the power of complex geometry, one likes to relate 
full CFT to CCFT. This is indeed possible. As far as the geometric 
aspects are concerned, the relationship is established by associating to $X$ 
its Schottky cover $\tC$, a double cover branched over the boundary $\partial X$
from which $X$ is recovered as the quotient by an anti-conformal involution. 
For $\partial X\eq\emptyset$ the double is just the total space of the 
orientation bundle. (For orientable boundaryless $X$ this bundle is trivial, 
i.e.\ the disconnected sum of two copies of $X$, with opposite orientations.) 
At the field theory level the relationship turns out to be quite a bit more
involved. But fortunately it can still be formulated in a model independent
manner. In particular, the quantities of basic physical interest are the
{\em \corfu s\/}, and these can be obtained as specific elements in 
spaces (\resp\ as sections in the bundles) of \cb s\ on $\tC$.
Which specific element in a block space constitutes a \corfu\
is determined by the following requirements:
\nxt
{}{\em Locality\/}: Correlation functions must be functions of the insertion 
points $p_s$
\\\mbox{}$\;\ \ $and (modulo the `Weyl anomaly') of the moduli of $X$.
\nxt
{}{\em Factorization\/}: Correlation functions on world sheets of different 
topology \\
\mbox{}$\;\ \ $must be compatible with the desingularization of singular curves $\tC$.

The factorization (or {\em gluing\/}, or {\em sewing\/}) conditions arise 
because one allows for singular curves $\tC$, possessing ordinary double points. 
When $\tC'$ is a partial desingularization of $\tC$ that resolves a 
double point $\pp\iN\tC$ in two points $p'{,}\,p''{\in}\,\hat C$,
chiral factorization provides a canonical isomorphism
$B_{\muv,\tC}\,{\cong}\,\bigoplus_{\nu\in\sus}\!B_{\muv,\nu,\nup;\tC'}$.  
This way one relates the chiral blocks at genus $g$ to the genus-zero blocks for 
  \be  \calh_{\mu_1}\,{\otimes}\,\calh_{\mu_2}\,{\otimes}\,\cdots
  \,{\otimes}\,\calh_{\mu_m} \otimes \Llb \mbox{$\bigoplus$}_{\lambda\in\sus}
  \calh_\lambda^{}{\otimes}\calh_\lambdap \Lrb^{\otimes g} \,. \end{equation}
The non-chiral factorization requires that the image of a \corfu\ under such an 
isomorphism is again a \corfu\ (for concrete formulas, see \cite{fffs3}). In 
physical terms, this amounts to a restriction on the allowed intermediate states
that contribute to singular limits of correlators, and is closely related to the
causality of dynamics and to the existence of operator product expansions.

The system of locality and factorization conditions is largely overdetermined.
Thus a priori the existence of solutions is questionable; but indeed, 
when the torus \parfu\ is of charge conjugation type and as long as there are 
no symmetry breaking boundary conditions (see below), existence can be proven 
by explicit construction via the relation with three-dimensional TFT 
\cite{fffs3}. There is also no guarantee of uniqueness. For closed orientable 
$X$ the solution is expected to be unique (and has been shown to be so for 
genus zero in large classes of models), once the modular invariant torus \parfu\
  \be  Z(\tau) = \sum_{\mu,\tilde\mu\in\sus} Z_{\mu,\tilde\mu}^{}\,\chii_\mu^{}
  (\tau)\,\overline{{\chii_{\tilde\mu}^{}(\tau)}_{}} \labl0z  \end{equation}
has been specified. In contrast, when one allows also for world sheets with
boundary, uniqueness requires the specification of
an allowed {\em \bc\/} on each boundary component.

\sect{Conformal boundary conditions}

The basic requirements in full CFT are the locality and factorization 
constraints. They allow e.g.\ to express all \corfu s on any closed orientable 
world sheet $X$ through the 3-point correlation functions on the sphere
(in string theory these determine the couplings of closed string \vop s) and
the 1-point correlation functions on the torus (including in particular
the torus \parfu\ \erf0z). When $X$ has a boundary\,%
 \footnote{~Such world sheets play e.g.\ a role in the analysis of
 dissipative quantum mechanics \cite{lesA}, of
 defects in systems of condensed matter physics, of percolation probabilities,
 and of string perturbation theory in the background of certain solitonic
 excitations, the so-called D-branes.}
and\,/\,or is non-orientable, 
there are further factorization conditions \cite{lewe3,fips,prss2,prss3}, which 
also account for the possible insertion of so-called {\em boundary fields\/}
\cite{card9}. (Boundary fields correspond to open string \vop s and can only be 
inserted on the boundary. The fields inserted in the interior of $X$ are called 
{\em bulk fields\/}.) Use of these constraints allows one to express all \corfu s
through the 3-point functions on the sphere, the torus \parfu, as well as:
\nxt
The 1-point \corfu s for bulk fields on $\mathbb{RP}^{\,2}$.
\nxt
The 3-point \corfu s for boundary fields on the disk.
\nxt
The \corfu s for one bulk and one boundary field on the disk.

\vskip.3em 
In physical terms, \corfu s\ are regarded as vacuum expectation values
of time-ordered products of field operators. Since the insertion point of a bulk 
field on $X$ has two pre-images in $\tC$, it accounts for two
marked points on $\tC$, to be labelled by two elements $\mu,\tilde\mu\iN\sus$.
In contrast, the insertion points for boundary fields lie 
on $\partial X$ and hence have a single pre-image in $\tC$,
so boundary fields carry a single chiral label $\mu\iN\sus$.
It is worth stressing that the boundaries of interest here are
{\em genuine, physical boundaries\/}. They must not be confused with 
the boundaries of small disks around field insertions that one sometimes cuts 
out from $X$ in order to specify local coordinates around the insertion points.

A major new issue are the {\em boundary conditions\/} that are to be specified
at such physical boundaries.\,%
 \footnote{~These, in turn, should not be mixed up
 with the periodic or twisted-periodic \bc s that are often studied
 in field theory. They specify the topological type of a
 bundle over the manifold and are not related to physical boundaries.} 
For investigating most aspects of \bc s one can concentrate on the
1-point functions $\langle\pho\mu\rangle$ for bulk fields on the disk $D$
(this is the basis of the so-called boundary state formalism).
As bulk fields carry two chiral labels $\mu,\tilde\mu\iN\sus$, these 
correspond to 2-point blocks $B_{\mu,\tilde\mu}$ on $\hat D\eq\Pe$ and,
due to ${\rm dim}(B_{\mu,\tilde\mu})\iN\{0,1\}$, are determined up to a scalar:
  \be  {\langle \pho\mu(v\ot\tilde v)\rangle}_{\!a}
  = \rca a\mu\vac B_{\mu,\tilde\mu}(v\ot\tilde v)  \,. \labl mm  \end{equation}
By factorization, the coefficients $\rca a\mu\vac$ appearing here are
interpreted as so-called {\em reflection coefficients\/}. These arise in 
the operator product that heuristically describes with which strength the
boundary vacuum field is excited when the bulk field $\pho\mu$ approaches
the boundary. One has $\rca a\mu\vac\eq0$ unless $\tilde\mu\eq\muP$; thus
we concentrate on the latter situation, and abbreviate 
${\rm R}^{a}_{\mu,\muP}\,{=:}\,{\rm R}^{a_{\phantom e}}_{\mu}$.

The new label $a$ in formula \erf mm indicates
that the 1-point function $\langle\pho\mu\rangle$ is not unique. Rather, it 
depends on the boundary condition to be attached to $\partial D$.
In other words, the label $a$ distinguishes
between distinct boundary conditions. A basic task is then to determine 
all conformally invariant \bc s that lead to \corfu s satisfying 
all consistency constraints. This task naturally splits into a chiral and a 
non-chiral part. At the chiral level, one wants to compute the `boundary blocks'
$B_{\mu,\tilde\mu}$ as solutions to the appropriate Ward identities.
The solution is known explicitly in many cases of interest, e.g.\ for
Neumann or Dirichlet conditions of free boson CFTs and for maximally
symmetric \bc s of \wzwm s. The resulting expressions are not particularly
useful, though. Rather, what is important are their normalization and their
factorization properties -- issues that can be studied model independently.

At the level of full CFT, the goal is to determine the reflection coefficients 
${\rm R}^a_\mu$. This can be achieved by solving a specific factorization 
identity, obtained by comparison of two singular limits of the 
\corfu\ $\langle\Pho\mu\Pho\nu{\rangle}_{\!a}$ for two bulk fields on $D$. 
An important observation is that this step is logically independent from the 
former. Thus one can classify \bc s even when the boundary blocks are not known
concretely. In fact this can be done in a model independent manner, only making 
use of the fact that the underlying CCFT obeys all chiral consistency conditions.

To remain within the framework of CFT, one wants the \bc s to preserve
the conformal symmetry, which means that the \cb s in the presence of the
boundary still satisfy the Ward identities associated to the Virasoro \alg.
But this requirement of conformal invariance turns out to be rather weak;
typically it is obeyed by infinitely many \bc s, so that a
classification is difficult. 
In contrast, a finite number of \bc s arises when one requires that they
preserve even a sub\alg\ $\chirb$ of \chir\ that is itself the VOA for some
rational CFT. In the sequel we assume to be in this situation and regard 
$\chirb$ to be given. The simplest case are 
the {\em bulk symmetry preserving\/} \bc s, for which $\chirb$ is all of $\chir$.

Let us be more specific. The appropriate labels of the bulk fields on the disk 
depend on the sub-VOA $\chirb$ that is prescribed to be preserved. But as will 
be seen below, they are not just given by the sector labels $\bar\mu$ of the
CCFT model that corresponds to $\chirb$; we denote the
label set by $\sush$ and its elements by $\Hat\mu,\Hat\nu,...$\,.
Independently of the choice of $\chirb$, factorization can be seen to imply 
  \be  \rch a\lambda\vac\, \rch a\mu\vac
  = \mbox{\Large$\sum$}^{}_{\Hat\nu\in\sush}\,
  \tNN\lambda\mu\nu\, \rch a\nu\vac \,,  \end{equation}
i.e.\ for every \bc\ $a$ the reflection coefficients $\rch a\mu\vac$ with
$\Hat\mu\iN\sush$
furnish a \onedim\ \irrep\ of some \alg, called the {\em \cla\/} and denoted
by \clAb. The structure constants $\tNN\lambda\mu\nu$ can be entirely 
expressed through CFT data that, manifestly, do not involve the boundary.

For $\chirb\eq\chir$ (bulk symmetry preserving \bc s), and with torus \parfu\ 
given by charge conjugation, i.e.\ $Z_{\lambda,\mu}^{}\,{=}\,\delta
_{\lambda,\muP}$, one finds that $\sush\eq\sus$ and
$\tNC\lambda\mu\nu\eq\N\lambda\mu\nu$, i.e.\ the
\cla\ \clA\ coincides with fusion \alg\ of the CCFT.
Thus in particular the basis elements of $\clA$ correspond to the sectors 
of the \chir-theory, the boundary labels $a$ are in the same set $\sus$,
and \clA\ is a semi-simple commutative associative \alg\ whose structure 
constants are expressible through a diagonalizing matrix $S$ as in \erf1v.

For $\chirb$ strictly contained in \chir\ -- {\em symmetry breaking\/} \bc s
-- the factorization arguments go through as well, so that one still deals with
\onedim\ \irrep s of a \cla. But the details are more involved. 
A systematic solution has been achieved \cite{fuSc112} for all cases where the
preserved sub-VOA \chirb\ is the fixed point \alg
  \be  \chirb_{\phantom |} = \chir^G  \end{equation} 
\wrt any finite abelian group $G$ of automorphisms of \chir,
or in other words, when \chirb\ is the VOA of an abelian orbifold
of the \chir-theory, as studied in section 5.
To determine the label set $\sush$ one considers the decomposition
  \be  \calh_\lambda^{} = \mbox{\Large$\bigoplus$}^{}_{\mub\in\bar\sus}\,
  V_\mub \otimes \calhb_\mub  \labl dc  \end{equation}
of simple \chir-modules as \chirb-modules.
Preservation of \chirb\ by the \bc\ means that the relevant \cb s respect the
Ward identities for \chirb, but not in general those for \chir, and hence the
sector labels $\mub$ of the \chirb-theory that appear in \erf dc -- corresponding
to untwisted orbifold sectors -- are to be
used. But in addition one and the same \chirb-module
$\bar\calh_\mub$ will typically give rise to different \cb s when it appears
in the decomposition of distinct \chir-modules
$\calh_\lambda$. One concludes that the labels $\Hat\mu$ are {\em pairs\/}
$(\mub{,}\psi_\mub)$ with a suitable degeneracy label $\psi_\mub$.

One must also determine the appropriate set of boundary labels $a$.
At the present state of affairs this still involves some heuristics.\,%
 \footnote{~This statement even applies to the symmetry preserving situation
 $\chirb\eq\chir$ as studied in \cite{card9} (except for some special classes of
 models where the structure constants $\tNC\lambda\mu\nu$ have been computed
 \cite{prss3} explicitly). In that case, however, many of the results can be
 made rigorous by employing the relation to TFT \cite{fffs3}. We expect that
 the methods of \cite{fffs3} can be extended so as to establish a proof also
 for more general \bc s. Work in this direction is in progress.}
The conclusions \cite{fuSc112} are that the boundary labels are $\calg$-orbits 
$[\rhob,\psu_\rhob]$ of pairs that consist of an (untwisted or twisted)
\chirb-sector label $\rhob$ and a character $\psu_\rhob$ of the central 
stabilizer $\calu_\rhob$, and that the degeneracy label $\psi_\mub$ of the labels 
in $\sush$ is a character of the (full) stabilizer $\cals_\mub$. Here $\calg$
is a simple current group that is necessarily present in the orbifold theory and 
which is naturally isomorphic to the character group $\Gs$ of the orbifold group 
$G$. One thus deals with two different deviations from the labels appearing in 
simple current extensions, cf.\ the text before \erf33. Namely, for $\hat\mu$ 
there is a character of the full rather than the central stabilizer and no 
$\calg$-orbit is to be formed, while $a$ has the same form as extended labels 
$\mu\xxt$, but the restriction to $T_{\J\mu}\eq T_\mu$ (that is, to untwisted 
orbifold sectors) is dropped. There are indications 
\cite{fuSc112,bppz2,fuSc14,scfu3} that \bc s violating this equality 
correspond to twisted or `solitonic' \rep s of the bulk symmetry \chir.

Further, the reflection coefficients are seen to be quotients
$\rch a\mu\vac\eq \tS_{\Hat\mu,a}/\tS_{\Hat\vac,a}$ with
  \be  \tS_{(\mub,\psi_\mub),{[\rhob,\psu_\rhob]}}
  = \frac{|G|^{\phantom|}_{}} {[\,|\cals_\mub|\,|\calu_\mub|\,
  |\cals_\rhob|\,|\calu_\rhob|\,{]}^{1/2}_{\phantom I}}
  \sum_{\J\in\cals_\mub\cap\calu_\rhob}
  \psi_\mub(\J)\, \psu_\rhob(\J)^*\, \SJ_{\mub,\rhob} \,.  \labl tS
  \end{equation}
Note that all input data appearing in this formula are purely chiral data.

Once these prescriptions for the labels and for the reflection coefficients 
are taken for granted, the \cla\ can be analysed rigorously. Among the
results \cite{fuSc112} are the following:
\Nxt
\clAb\ is a unital semi-simple commutative associative \alg.
\nxt
The unit element is $\Hat\vac\eq\vacb$, the vacuum sector of the $G$-orbifold.
\nxt
The structure constants of \clAb\ obey the Verlinde-like formula
  \be  \tNl\lambda\mu\nu =
  \sum_{a}  \frac{ \tS_{\Hat\lambda,a}\, \tS_{\Hat\mu,a}\,
  \tS_{\Hat\nu,a}} {\tS_{\hat\vac,a}} \,.  \end{equation}
(To raise the third index, one contracts with the inverse of
$\hat C_{\hat\lambda,\hat\mu}\eq\tNl\lambda\mu\vac$.)
\nxt
\clAb\ contains the fusion rule \alg\ ${\cal A}(\chir)$ as a sub\alg.
\nxt
More specifically, \clAb\ is a direct sum of ideals $\calc^\sigma(\chir)$
with $\sigma\iN G$. Each of them corresponds to \bc s of a definite
{\em \auto\ type\/} $\sigma$ (and depends only on the automorphism $\sigma$, 
but not on the specific orbifold group $G\,{\ni}\,\sigma$ considered).
The latter is characterized by the fact that in the
presence of such a \bc\ the 1-point functions of bulk fields 
can be expressed entirely through {\em twisted boundary blocks\/}
$B^\sigma_{\mu,\muP}$ which satisfy $\sigma$-twisted Ward identities.
Schematically, 
  \be  B^\sigma_{\mu,\muP} \circ
  \llb Y_n \oT\bfe - (-1)^{\Delta_Y}_{}\,\bfe\oT \sigma(Y_{\!-n})\lrb
  = 0  \end{equation}
for every field $Y$, of conformal weight $\Delta_Y$, in the VOA \chir.
Twisted boundary blocks are related to ordinary boundary blocks 
$B_{\mu,\mup}$ as
  \be  B_{\mu,\mup}^\sigma = B_{\mu,\mup}^{} \circ (\Tau_\sigma \oT \bfe\,) 
  \end{equation}
with $\Tau_\sigma$ the twisted intertwiner associated to $\sigma$.

\smallskip

Let us add a few further comments:
\Nxt
Symmetry breaking \bc s play e.g.\ an important role for space-time 
supersymmetry breaking in superstring theory. In that context the preserved 
sub\alg\ \chirb\ must contain the $N\eq1$ world sheet superconformal 
symmetry, which is a  gauge symmetry of (perturbative) superstring theory.
On the other hand, in condensed matter applications
one may wish to study \bc s that even break the conformal symmetry; these
are not covered by the CFT constructions reported above.
\nxt
The \bc s preserving $\chirg\,{\subset}\,\chir$ include of course all those
which even preserve some intermediate \alg\ $\chir'$ satisfying
$\chirg\,{\subset}\,\chir'\,{\subset}\,\chir$. By
the Galois theory of \voa s \cite{doMa5}, there is a bijection between these 
intermediate \alg s and the set of subgroups of $G$.
\nxt 
That every \bc\ possesses an \auto\ type is no longer true in general when 
$\chirb\nE\chirg$. This property is only to be expected when $\chirb$ is an 
orbifold subalgebra, and it has only been proven for abelian orbifold groups 
$G$. In string theory, \bc s without \auto\ type correspond to so-called
non-BPS D-branes (for a review of the latter see e.g.\ \cite{Leru}).
\nxt
The one-point functions on $\mathbb{RP}^{\,2}$ have also been studied
in the literature, see e.g.\ \cite{fips,prad,huss2}. 
Recently, closed formulas for the one-point functions on the disk and on
$\mathbb{RP}^{\,2}$ have been found \cite{fhssw} for 
the case of arbitrary torus \parfu s that are related to simple currents.
\nxt 
There is evidence \cite{scfu3} that for arbitrary (rational) $\chirb$
and arbitrary torus partition function the set of \bc s preserving
$\chirb$ comes equipped with the structure of a spherical tensor category, 
called the `boundary category', and that this category can be entirely
constructed from the modular tensor category associated to $\chirb$.
The structures arising in this context are familiar from subfactor 
theory \cite{boev5} and suggest an interpretation of the symmetry
breaking \bc s in terms of suitable `solitonic' \rep s of the bulk chiral 
\alg\ \chir. (When \chirb\ is an orbifold sub\alg\ $\chirg$ of \chir,
then these solitonic \rep s are nothing but the usual twisted sector \rep s
\cite{dolm3} of the orbifold model.)
\\
The same arguments as for symmetry preserving \bc s show that the operator
product of two boundary fields corresponds to the $6j$-symbols of the
boundary category. The fusion rules of the boundary category coincide
with the annulus coefficients, i.e.\ the coefficients in an expansion of
the \parfu\ on the annulus \wrtt characters of the solitonic \rep s.

\smallskip

Finally let us illustrate our results by exposing the 
boundary conditions of \wzwm s that preserve at least the fixed point \alg\
of \chir\ under a $\zet_2$ orbifold group $G\eq\{{\rm id},\sigma\}$. 
(Other WZW-orbifolds are rather more complicated, but for the analysis of \bc s
fortunately only partial results about the orbifolds are needed.)
Using the information about the
corresponding WZW-orbifolds collected in section 5, one finds that the 
elements of $\sush$ can be suggestively written by using sector labels of
the full bulk symmetry, according to
  \be  \lambda \;\ {\rm and}\;\ \sigma^*\!\lambda
  \quad {\rm for}\ \sigma^*\!\lambda\nE\lambda \,, \quad\ \mbox{and}\quad
  (\lambda,\psi) \ (\psi\in\{\pm1\})
  \;\ {\rm for}\ \sigma^*\!\lambda\eq\lambda  \,.  \end{equation}
Similarly, the labels for the \bc s are
  $$   \bearl
  \mu \quad \mbox{for length-2 orbits}\ \{(\mu,1,0),(\mu,-1,0)\} 
  \ \mbox{of untwisted orbifold sectors} \,,
  \\{}\\[-.8em]
  \mu\;\ {\rm and}\;\ \sigma^*\!\mu \quad \mbox{for $\calg$-fixed points}\
  (\mu,0,0) \,,
  \\{}\\[-.8em]
  {\mud} \quad \mbox{for length-2 orbits}\ \{(\mud,1,1),(\mud,-1,1)\}
  \ \mbox{of twisted orbifold sectors}\,.  \eear $$ 
In this notation, the entries of the diagonalizing matrix $\tS$ read
  \be \!\!\! \bearll
  \left.\bearl  {\tS_{(\lambda,\psi),\mu} = S_{\lambda,\mu}} \\{}\\[-.7em]
  {\tS_{(\lambda,\psi),\mud} = \psi\eta_\lambda^{-1}\,\Se_{\lambda,\mud}}
  \eear\!\!\right\}\,\mbox{for $\sigma^*\!\lambda\eq\lambda$} \,, &
  \left.\bearl  {\tS_{\lambda,\mu} = S_{\lambda,\mu}} \\{}\\[-.7em]
  \tS_{\lambda,\mud} = 0
  \eear\!\!\right\}\,\mbox{for $\sigma^*\!\lambda\nE\lambda$} \,. \eear
  \end{equation}
For the structure constants of the \cla\ one finds e.g.\
  \be \TN_{\lambda_1,\lambda_2}^{\;\ \ \ \lambda_3}
  \eq \N{\lambda_1}{\lambda_2}{\;\lambda_3} \,, \quad\ 
  \TN_{\lambda_1,\lambda_2}^{;\ \ (\lambda_3,\psi_3)}
  \eq \N{\lambda_1}{\lambda_2}{\;\lambda_3} \,, \quad\ 
  \TN_{\lambda_1,(\lambda_2,\psi_2)}^{\ \ \ \ (\lambda_3,\psi_3)}
  \eq \N{\lambda_1}{\lambda_2}{\;\lambda_3} \,,  \end{equation}
as well as
  \be  \bearll
  \TNl{\lambda_1}{\psi_1}{\lambda_2}{\psi_2}{\lambda_3}{\psi_3} \!\!\!
  &= {\rm N}_{\lambda_1,\lambda_2,\lambda_3} \,{+}\, \Frac{\psi_1\psi_2\psi_3}
  {\eta^{}_{\lambda_1} \eta^{}_{\lambda_2}\eta^{}_{\lambda_3}}\,
  {\rm N}^{\sss(0)}_{\lambda_1,\lambda_2,\lambda_3}
  \\{}\\[-.82em]  &= {\rm N}^{\cal O}
  _{(\lambda_1,\psi_1,0),(\lambda_2,\psi_2,0),(\lambda_3,\psi_3,0)} \,,
  \eear \end{equation}
with ${\rm N}^{\sss(0)}_{\lambda_1,\lambda_2,\lambda_3}\,{:=}\,\sum_\mud 
\Se_{\lambda_1,\mud}\Se_{\lambda_2,\mud}\Se_{\lambda_3,\mud}/\Se_{\vac,\mud}$.  


 \newcommand\wb{\,\linebreak[0]} \def\wB {$\,$\wb}
 \newcommand\Bi[2]    {\bibitem[#2]{#1}}
 \renewcommand\J[5]   {{\em #5}, {#1} {\bf #2} ({#3}), {#4} }
 \newcommand\JK[5]    {{\em #5}, {#1} {\bf #2} ({#3}) }
 \newcommand\Prep[2]  {{\em #2}, pre\-print {#1}} 
 \newcommand\PreP[2]  {{\em #2}, pre\-print, {#1}} 
 \newcommand\PRep[2]  {{\em #2}, {#1}} 
 \newcommand\BOOK[4]  {{\em #1\/} ({#2}, {#3} {#4})}
 \newcommand\inBO[7]  {{\em #7}, in:\ {\em #1}, {#2}\ ({#3}, {#4} {#5}), p.\ {#6}}
 \newcommand\Erra[3]  {\,[{\em ibid.}\ {#1} ({#2}) {#3}, {\em Erratum}]}
 \def\jf    {J.\ Fuchs}
 \def\adma  {Adv.\wb Math.}
 \def\anop  {Ann.\wb Phys.}
 \def\aspm  {Adv.\wb Stu\-dies\wB in\wB Pure\wB Math.}
 \def\atmp  {Adv.\wb Theor.\wb Math.\wb Phys.}
 \def\coma  {Con\-temp.\wb Math.}
 \def\comp  {Com\-mun.\wb Math.\wb Phys.}
 \def\cpma  {Com\-pos.\wb Math.}
 \def\duke  {Duke\wB Math.\wb J.}
 \def\gafa  {Geom.\wB and\wB Funct.\wb Anal.}
 \def\ijmp  {Int.\wb J.\wb Mod.\wb Phys.\ A}
 \def\jgap  {J.\wb Geom.\wB and\wB Phys.}
 \def\joag  {J.\wB Al\-ge\-bra\-ic\wB Geom.}
 \def\joal  {J.\wB Al\-ge\-bra}
 \def\jomp  {J.\wb Math.\wb Phys.}
 \def\lemp  {Lett.\wb Math.\wb Phys.}
 \def\maan  {Math.\wb Annal.}
 \def\maze  {Math.\wb Zeitschr.}
 \def\nuci  {Nuovo\wB Cim.}
 \def\nupb  {Nucl.\wb Phys.\ B}
 \def\phep  {Proc.\wb HEP}
 \def\phlb  {Phys.\wb Lett.\ B}
 \def\phrl  {Phys.\wb Rev.\wb Lett.}
 \def\rims  {Publ.\wB RIMS}
 \def\sebo  {S\'emi\-naire\wB Bour\-baki}
 \def\slnp  {Sprin\-ger\wB Lecture\wB Notes\wB in\wB Physics}
 \newcommand\mbop[2] {\inBO{The Mathematical Beauty of Physics}
            {J.M.\ Drouffe and J.-B.\ Zuber, eds.} \WS\Si{1997} {{#1}}{{#2}} }
 \newcommand\geap[2] {\inBO{Physics and Geometry} {J.E.\ Andersen, H.\
            Pedersen, and A.\ Swann, eds.} \MD\NY{1997} {{#1}}{{#2}} }
   \def\AMS    {{American Mathematical Society}}
   \def\AW     {{Addi\-son\hy Wes\-ley}}
   \def\Be     {{Berlin}}
   \def\Ca     {{Cambridge}}
   \def\CUP    {{Cambridge University Press}}
   \def\MD     {{Marcel Dekker}}
   \def\NY     {{New York}}
   \def\OUP    {{Oxford University Press}}
   \def\PR     {{Providence}}
   \def\SV     {{Sprin\-ger Ver\-lag}}
   \def\WS     {{World Scientific}}
   \def\Si     {{Singapore}}

\end{document}